\documentclass[11pt]{amsart}
\usepackage{xypic}
\usepackage{amssymb,amsmath,amsgen,amsxtra,amsfonts} 
\usepackage{dsfont,times}

\begin{document}
%
%   D e f i n i t i o n s
%
% gibt's die folgenden Definition schon?
%\newcommand{\qed}{\mbox{}\hfill$\Box$}
%\newcommand{\mod}{{\:\mbox{mod}\:}}
%
\theoremstyle{definition}
\newtheorem{Definition}{Definition}[section]
\newtheorem{Convention}{Definition}[section]
\newtheorem{Construction}{Construction}[section]
\newtheorem{Example}[Definition]{Example}
\newtheorem{Examples}[Definition]{Examples}
\newtheorem{Remark}[Definition]{Remark}
\newtheorem{Remarks}[Definition]{Remarks}
\newtheorem{Caution}[Definition]{Caution}
\newtheorem{Conjecture}[Definition]{Conjecture}
\newtheorem{Question}[Definition]{Question}
\newtheorem{Questions}[Definition]{Questions}
\theoremstyle{plain}
\newtheorem{Theorem}[Definition]{Theorem}
\newtheorem*{Theoremx}{Theorem}
\newtheorem*{Definitionx}{Definition}
\newtheorem{Proposition}[Definition]{Proposition}
\newtheorem{Lemma}[Definition]{Lemma}
\newtheorem{Corollary}[Definition]{Corollary}
\newtheorem{Fact}[Definition]{Fact}
\newtheorem{Facts}[Definition]{Facts}
\newtheoremstyle{voiditstyle}{3pt}{3pt}{\itshape}{\parindent}%
{\bfseries}{.}{ }{\thmnote{#3}}%
\theoremstyle{voiditstyle}
\newtheorem*{VoidItalic}{}
\newtheoremstyle{voidromstyle}{3pt}{3pt}{\rm}{\parindent}%
{\bfseries}{.}{ }{\thmnote{#3}}%
\theoremstyle{voidromstyle}
\newtheorem*{VoidRoman}{}

% abgeschrieben aus The LaTeX Companion, 2nd edition,
% von Mittelback & Goossens
%
\newcommand{\prf}{\par\noindent{\sc Proof.}\quad}
\newcommand{\blowup}{\rule[-3mm]{0mm}{0mm}}
\newcommand{\cal}{\mathcal}
\newcommand{\Aff}{{\mathds{A}}}
\newcommand{\BB}{{\mathds{B}}}
\newcommand{\CC}{{\mathds{C}}}
\newcommand{\FF}{{\mathds{F}}}
\newcommand{\GG}{{\mathds{G}}}
\newcommand{\HH}{{\mathds{H}}}
\newcommand{\NN}{{\mathds{N}}}
\newcommand{\ZZ}{{\mathds{Z}}}
\newcommand{\PP}{{\mathds{P}}}
\newcommand{\QQ}{{\mathds{Q}}}
\newcommand{\RR}{{\mathds{R}}}
\newcommand{\Liea}{{\mathfrak a}}
\newcommand{\Lieb}{{\mathfrak b}}
\newcommand{\Lieg}{{\mathfrak g}}
\newcommand{\Liem}{{\mathfrak m}}
\newcommand{\ideala}{{\mathfrak a}}
\newcommand{\idealb}{{\mathfrak b}}
\newcommand{\idealg}{{\mathfrak g}}
\newcommand{\idealm}{{\mathfrak m}}
\newcommand{\idealp}{{\mathfrak p}}
\newcommand{\idealq}{{\mathfrak q}}
\newcommand{\idealI}{{\cal I}}
\newcommand{\lin}{\sim}
\newcommand{\num}{\equiv}
\newcommand{\dual}{\ast}
\newcommand{\iso}{\cong}
\newcommand{\homeo}{\approx}
\newcommand{\mm}{{\mathfrak m}}
\newcommand{\pp}{{\mathfrak p}}
\newcommand{\qq}{{\mathfrak q}}
\newcommand{\rr}{{\mathfrak r}}
\newcommand{\pP}{{\mathfrak P}}
\newcommand{\qQ}{{\mathfrak Q}}
\newcommand{\rR}{{\mathfrak R}}
%
%  evtl. auch \"uber \mathbb oder \Bbb
%
\newcommand{\dq}{{``}}
\newcommand{\OO}{{\cal O}}
\newcommand{\into}{{\hookrightarrow}}
\newcommand{\onto}{{\twoheadrightarrow}}
\newcommand{\Spec}{{\rm Spec}\:}
\newcommand{\uSpec}{{\bf Spec}\:}
\newcommand{\Proj}{{\rm Proj}\:}
\newcommand{\Pic}{{\rm Pic }}
\newcommand{\Br}{{\rm Br}}
\newcommand{\fBr}{{\widehat{\rm Br}}}
\newcommand{\NS}{{\rm NS}}
\newcommand{\chit}{\chi_{\rm top}}
\newcommand{\KanDiv}{{\cal K}}
\newcommand{\perdef}{{\stackrel{{\rm def}}{=}}}
%plan abschaffung du
\newcommand{\Cycl}[1]{{\ZZ/{#1}\ZZ}}
\newcommand{\Sym}{{\mathfrak S}}
\newcommand{\Xcan}{X_{{\rm can}}}
\newcommand{\ab}{{\rm ab}}
\newcommand{\Aut}{{\rm Aut}}
\newcommand{\Hom}{{\rm Hom}}
\newcommand{\shHom}{{\cal Hom}}
\newcommand{\Supp}{{\rm Supp}}
\newcommand{\ord}{{\rm ord}}
\newcommand{\divisor}{{\rm div}}
\newcommand{\Alb}{{\rm Alb}}
\newcommand{\Jac}{{\rm Jac}}
\newcommand{\piet}{{\pi_1^{\rm \acute{e}t}}}
\newcommand{\Het}[1]{{H_{\rm \acute{e}t}^{{#1}}}}
\newcommand{\Hcris}[1]{{H_{\rm cris}^{{#1}}}}
\newcommand{\HdR}[1]{{H_{\rm dR}^{{#1}}}}
\newcommand{\hdR}[1]{{h_{\rm dR}^{{#1}}}}
\newcommand{\Mev}{{\cal M}^{\rm even}}
\newcommand{\Modd}{{\cal M}^{\rm odd}}
\newcommand{\Mevi}{{\cal M}^{\rm even, insep}}
\newcommand{\Moddi}{{\cal M}^{\rm odd, insep}}
\newcommand{\defin}[1]{{\bf #1}}
\title[Horikawa Surfaces in Positive Characteristic]{The Canonical Map and Horikawa Surfaces in Positive Characteristic}

%\keywords{Surfaces of general type, inseparable morphism, Enriques--Horikawa surfaces}

\author[Christian Liedtke]{Christian Liedtke}
\address{Department of Mathematics, Stanford University, 450 Serra Mall, Stanford, CA 94305, USA}
\curraddr{}
\email{liedtke@math.stanford.edu}

\date{December 20, 2011}

\begin{abstract}
We extend fundamental inequalities related to the canonical map of surfaces of 
general type 
%due to Beauville, Castelnuovo, Horikawa and Noether 
to positive characteristic.
Next, we classify surfaces on the Noether lines, i.e., even and odd Horikawa
surfaces, in positive characteristic.
We describe their moduli spaces and the subspaces formed by surfaces whose canonical
maps are inseparable.
Moreover, we compute their Betti-, deRham- and crystalline cohomology.
Finally, we prove lifting to characteristic zero and show that the moduli
spaces are topologically flat over the integers. 
\end{abstract}
\setcounter{tocdepth}{1}
\maketitle
%\tableofcontents
%
\section*{Introduction}

A detailed classification of surfaces of general type, even over the complex numbers,
still looks hopeless.
A more modest aim is to prove (in-)equalities holding between
their fundamental invariants.
Next, one would like to classify at least those surfaces that are extremal with respect
to these inequalities.
Even over the complex numbers, this is still far from complete and nearly nothing
is known over fields of positive characteristic.

One of these fundamental inequalities for minimal surfaces of
general type is
 \begin{center}
 \begin{tabular}{clc}
   \blowup %{\rm ($\ast$)}
   & {\rm Noether's inequality} & $ K_X^2\,\geq\,2p_g(X)-4 $ 
 \end{tabular}
 \end{center}
see \cite{noether} and \cite[Chapter VII]{bhpv} 
for the complex case and \cite{lie2} for
arbitrary characteristic.
It generalizes the formula $\deg K_C=2g-2$ for curves.

Depending on the image of the canonical map, there are sharper inequalities
by Beauville \cite{bea} and Castelnuovo \cite{castelnuovo}.
As our first result we extend these to arbitrary characteristic

\begin{Theoremx}
 Let $X$ be a minimal surface of general type in arbitrary characteristic.
 In case the canonical image is
 \begin{center}
 \begin{tabular}{llc}
   \blowup %{\rm (1)} & 
     a curve, & {\em Beauville's inequality} & $K_X^2\,\geq\,3p_g(X)-6 $ \\
   \blowup %{\rm (2)} & 
     birational to $X$, & {\em Castelnuovo's inequality} & $K_X^2\,\geq\,3p_g(X)-7$
 \end{tabular}
 \end{center}
 holds true.
\end{Theoremx}

Not too surprisingly, the main difficulty lies in establishing these
inequalities in small characteristics.
As in the complex case, the proofs give a
structure result for minimal surfaces of general type with 
$K^2\leq3p_g-7$.
We refer to Section \ref{sec:canonical map} for the precise statement.

Coming back to Noether's inequality, %($\ast$) 
we recall that a minimal surface of general type is called an
\begin{center}
 \begin{tabular}{ll}
  \blowup {\em even Horikawa surface} & if \quad $K_X^2\,=\,2p_g(X)-4$, \quad and\\
  \blowup {\em odd Horikawa surface}  & if \quad $K_X^2\,=\,2p_g(X)-3$
 \end{tabular}
\end{center}
holds.
Both classes together are referred to as 
{\em surfaces on the Noether lines}.
Over the complex numbers, they have been roughly classified 
by Enriques \cite{enr} and more thoroughly by 
Horikawa \cite{horikawa}, \cite{horikawa2}.
A characteristic-free classification of even
Horikawa surfaces has already been obtained in \cite{lie2}.
Here, we reprove this result and extend it also to
odd Horikawa surfaces.

\begin{Theoremx}
 Let $X$ be a surface on the Noether lines in arbitrary characteristic.
 \begin{enumerate}
  \item If $X$ is an even Horikawa surface then the canonical linear
   system is base point free.
  \item If $X$ is an odd Horikawa surface with $p_g\geq5$ then the canonical
   linear system has a unique base point, which can be resolved after one 
   blow-up. 
 \end{enumerate}
 In both cases the canonical map is generically finite of degree $2$.
 Its image is a possibly singular rational surface of degree $p_g-2$ in 
 $\PP^{p_g-1}$.
\end{Theoremx}

Using del~Pezzo's results on surfaces of degree $n-1$ in $\PP^n$ it
is not difficult to work out the detailed classification, 
see Section \ref{sec:noether}.
Next, we prove existence:

\begin{Theoremx}
  All cases of the previous structure result
  exist in arbitrary characteristic.
  Moreover, in characteristic $2$ all cases exist with inseparable 
  as well as separable canonical maps.
\end{Theoremx}

However, surfaces with inseparable canonical maps exist only in
characteristic $2$ and even there they are quite rare.
More precisely,

\begin{Theoremx}
 Surfaces with inseparable canonical maps form a proper closed subset of
 the moduli space of even Horikawa surfaces with fixed canonical image.
 Both moduli spaces are irreducible and unirational. 
\end{Theoremx}

In fact, we determine the codimension in Section \ref{sec:moduli}.
The canonical model of a generic even Horikawa surface with inseparable 
canonical map has usually lots of $A_1$-singularities,
whereas the canonical model of a generic even Horikawa surface tends to 
be smooth.
We refer to Section \ref{sec:moduli} for precise and quantitative
statements.

Odd Horikawa surfaces with $p_g\leq4$ fall into more
classes and we refer to Section \ref{sec:small} for a detailed and
characteristic-free classification.
For example, we find surfaces with $p_g=K^2=3$ and inseparable
canonical maps in characteristic $3$.

Using our explicit classification we compute the basic invariants

\begin{Theoremx}
 Let $X$ be a surface on the Noether lines.
 Then
 \begin{enumerate}
  \item it is algebraically simply connected,
  \item its Picard scheme is reduced,
  \item the Fr\"olicher spectral sequence 
%from Hodge- to deRham-cohomology
    degenerates at $E_1$-level,
  \item its crystalline cohomology is torsion-free, and
  \item the Betti- and deRham-numbers coincide.
 \end{enumerate} 
\end{Theoremx}

On the other hand, we give examples in arbitrary large characteristics
where the slope spectral sequence from Hodge--Witt to crystalline
cohomology does not degenerate at $E_1$-level.

Finally, we prove lifting:

\begin{Theoremx}
 Let $X$ be a surface on the Noether lines over 
 $k=\overline{k}$.
 Then
 \begin{enumerate}
  \item its canonical model lifts over the Witt ring $W(k)$, and 
  \item $X$ lifts via an algebraic
    space over a possibly ramified extension of $W(k)$.
 \end{enumerate}
\end{Theoremx}

As a direct consequence we obtain

\begin{Theoremx}
 The moduli spaces of surfaces on the Noether lines are
 topologically flat over $\Spec\ZZ$.
\end{Theoremx}

Unfortunately, in order to prove flatness and to determine the
scheme structure of these moduli spaces, we have to understand
contributions of singularities to the deformation theory 
of canonical models better.
This will be dealt with in future work.
\medskip 

The article is organized as follows:

In Section \ref{sec:canonical map} we extend
Beauville's results circling
around Castelnuovo's inequality \cite[Section 5]{bea}
to arbitrary characteristic.
In particular, we obtain the structure result for surfaces
with $K^2\leq3p_g-7$.

In Section \ref{sec:alphaltorsors} we introduce an easy method
to construct inseparable morphisms of degree $p$
using generic $\alpha_{\cal L}$-torsors.
We use Cartier's canonical connection 
to describe their singularities
and prove Bertini-type theorems on generic singularities.
This section is independent from the rest of the article
and will be useful for the construction of surfaces
in many other contexts.

From Section \ref{sec:noether} on we discuss surfaces on the Noether
lines. 
In Section \ref{sec:noether} their basic structure result is proven
(assuming $p_g\geq5$ for odd Horikawa surfaces) and 
existence is shown in arbitrary characteristic.
In characteristic $2$ we show that all types of
surfaces exist with separable as well as inseparable canonical maps.

In Section \ref{sec:moduli} we give a characteristic-free 
description of the moduli space of even Horikawa surfaces with 
fixed canonical image.
Moreover, we describe the proper subset formed by surfaces
with inseparable canonical map.
Finally, we determine the singularities on generic
canonical models in these moduli spaces.

In Section \ref{sec:small} we classify the remaining surfaces
on the Noether lines,
namely odd Horikawa surfaces with $p_g\leq4$ and 
establish their existence.
%in arbitrary characteristic.
%This section is separate since extra cases show up for these
%small values of $p_g$.

In Section \ref{sec:hodge} we show that the Fr\"olicher spectral
sequence 
%from Hodge- to deRham-cohomology 
for surfaces on the Noether lines degenerates at $E_1$.
Next, we compute their Betti-numbers, deRham-numbers and deduce
that their crystalline cohomology groups are torsion free.
Finally, we give examples in arbitrary large characteristic
where the slope spectral sequence does not degenerate at $E_1$.

In Section \ref{sec:lifting} we prove lifting of the canonical models 
over the Witt ring.
This immediately implies that their moduli spaces
are topologically flat over $\Spec\ZZ$.

Finally, Appendix \ref{sec:appendix} contains a couple of results
on double covers, which are somewhat scattered over the literature.
Also, we extend results on canonical resolution of Du~Val singularities
to characteristic $2$.

\begin{VoidRoman}[Acknowledgements]
  I thank Fabrizio Catanese for a lively discussion at the 
  University of Warwick, which
  led to the results in Section \ref{sec:moduli}.
  Also, I thank Torsten Ekedahl and the referee for comments.
  Finally, I gratefully acknowledge funding from DFG under 
  research grant LI 1906/1-1 and thank the department 
  of mathematics  at Stanford university for kind hospitality.
\end{VoidRoman}

\section{The canonical map}
\label{sec:canonical map}

In this section we establish several classical results on the canonical map of a
surface of general type in positive characteristic.
Basically, the classical proofs also work in positive characteristic, although 
a little care is needed since generic fibers of morphisms may no longer be smooth
%(Bertini's theorem is no longer valid),
and the vanishing results of Kodaira, Mumford and Ramanujam are not at
disposal.

\subsection*{First case: the canonical image is one-dimensional}
Our first result is due to Beauville \cite[Lemme 5.3]{bea} over the complex numbers.
The main difficulty extending this result to positive characteristic are
inseparability possibilities in characteristic $2$, on which we comment during
our proof.

\begin{Theorem}
\label{pencil}
  Let $X$ be a minimal surface of general type.
  If the canonical linear system of $X$ defines a rational map with one-dimensional image,
  then {\em Beauville's inequality}
  $$
     K_X^2\,\geq\,3p_g\,-\,6
  $$
  holds true.
\end{Theorem}

\prf
Let $\widetilde{X}\to X$ be the resolution of indeterminacy
and $p:\widetilde{X}\to B$ be the Stein factorization of the canonical map.
The generic fiber is a connected and geometrically integral curve $F$.
% Badescu's book, Theorem 7.1.
We denote its arithmetic genus by $p_a(F)$.

If the canonical map has basepoints, then $B$ is rational.
Writing $|K_X|=Z+|aF|$, % where $F$ is a general irreducible fiber,
Riemann--Roch on $B$ yields $a=p_g-1$.
Moreover, we calculate
$$
K_X^2\,\geq\,a\,K_XF\,\geq\,a^2\,F^2\,=\,(aF^2)\,(p_g-1)
$$
We are done if $aF^2\geq3$ and the result is trivial for $a=1$.
In the remaining case $a=2$, $F^2=1$ we find $K_X^2\geq4=3p_g-5$.
This establishes the desired inequality.

Thus, we may assume that the canonical linear system defines a morphism.
We write it as $|K_X|=Z+|aF|$, where $Z$ denotes
the fixed part.
Hence there exists a divisor $D_a$ of degree $a$ on the base $B$ 
such that $K_X\lin Z+p^\ast D_a$.
By Riemann--Roch on $B$ we obtain $a\geq p_g + p_a(B) - 1$.

In case $p_a(F)\geq3$ we are done since then
$$
K_X^2\,\geq\,a\, K_XF\,\geq\,4a\,\geq\, 4(p_g-1)\,.
$$

Hence we may assume $p_a(F)=2$. 
We write the canonical system as $|K_X|=V+H+|aF|$, where $V+H$ denotes the fixed
part of the canonical system decomposed into a vertical component $V$ and
a horizontal component $H$.
From $K_XF=2$ we obtain $HF=2$.

In case $H$ is reduced, it is an irreducible curve
or the the sum of two irreducible curves.
Hence we obtain
$$
\begin{array}{lcll}
  HK_X\,+\,H^2 &\geq& -4 &\mbox{ by the adjunction formula}\\
  HK_X-H^2&=&HV+2a\geq2a
\end{array}
$$
Thus, $HK_X\geq a-2$ and we obtain
$$
K_X^2\,\geq\,HK_X\,+\,2a\,\geq\,3a-2\,\geq\,3p_g-5
$$ 
It remains to treat the case $H=2C$, where $C$ is an irreducible curve that
defines a section of the fibration.

So far we followed Beauville's proof \cite[Lemme 5.3]{bea}.
However, the generic fiber may not be smooth, which can happen only in
characteristic $p\leq5$ by Tate's theorem on
genus change in inseparable field extensions and $p_a(F)=2$.
Moreover, in characteristic $2$ it could happen that the map of degree $2$
onto a ruled surface that he constructs becomes inseparable.

Instead, we consider $p_\ast\omega_X$ on $B$.
For every $b\in B$ the fiber $F_b=p^{-1}(b)$
has a reduced component (the curve $C$ defines a section),
and hence $h^0(F_b,\OO_{F_b})=1$.
Since the Euler characteristic is constant in flat families
we obtain $h^0(F_b,\omega_{F_b})=2$ for all $b\in B$.
The adjunction formula yields
$$
\omega_X|_{F_b} \,\iso\,\left(\omega_X\otimes\OO_X(F_b)\right)|_{F_b}\,\iso\,\omega_{F_b}\,.
$$
The global sections being $2$-dimensional in {\em every} fiber,
Grauert's theorem 
\cite[Corollary III.12.9]{hart}
implies that $p_\ast\omega_X$ is locally free of rank $2$ 
on $B$ and that for all $b\in B$ the natural map
$$
p_\ast\omega_X\,\otimes_{\OO_B}\,k(b) \,\to\, H^0(F_b,\,\omega_X|_{F_b})
$$
is an isomorphism.

Hence the natural morphism $p^\ast p_\ast\omega_X\to\omega_X$ defines a rational
map $X\dashrightarrow S:=\PP(p_\ast\omega_X)$, where $S$ is a ruled surface
over $B$.
On the generic fiber over $B$ it induces the
canonical map, which is finite of degree $2$ onto $\PP^1_{k(B)}$
(the generic fiber is a regular curve with $p_a=2$).
Thus, $X\dashrightarrow S$ is a generically finite, possibly only 
rational map of degree $2$.

For every $b\in B$ the fiber $F_b$ is a Gorenstein curve with 
$h^0(\omega_{F_b})=2$.
Intersecting $C$ with $F_b$ we obtain a Cartier divisor $C_b$ of degree $1$
on $F_b$.
Furthermore, we have $h^0(\omega_{F_b}(-C_b))\leq1$, which follows from
Riemann--Roch in case $h^1$ of this sheaf is zero and from Clifford's
inequality in case $h^1$ does not vanish.
We stress that this is true for every fiber since the theorems of
Riemann--Roch and Clifford hold for Gorenstein curves.
% Clifford fuer reduzible Kurven ?
% h^0=h^1 da \chi=0 
% problematisch wird's nur, wenn h^0 \geq 2
% Kohomologie ist supported auf der reduziblen Komponente und
% der Clifford-Beweis geht klar.

For a suitable blow-up $\epsilon:\widehat{X}\to X$, the rational
map extends to a morphism $\pi:\widehat{X}\to S$.
From the previous discussion we see that for every $b\in B$, the Cartier 
divisor $C_b$ on $F_b$ is not a base point for $\omega_{F_b}$.
In particular, we do not have to blow up points on $C$ in the resolution
$\epsilon$ of indeterminacy and it follows that the total transform
$\epsilon^\ast(C)=\widehat{C}$ is an irreducible curve on $\widehat{X}$.

Now, we follow again \cite[Lemme 5.3]{bea}.
From the case-by-case-analysis above we still have
$$
K_{\widehat{X}}\,=\,\widehat{V}\,+\,2\widehat{C}\,+\,\widehat{p}^\ast D_a,
$$
where $\widehat{V}$ is a vertical divisor for $\widehat{p}$, and where
$\widehat{V}+\widehat{C}$ is the fixed part of $|K_{\widehat{X}}|$.
On the other hand, we know that
$$
|K_{\widehat{X}}|\,=\,\pi^\ast| K_S+\delta | 
$$
for some divisor $\delta$ on $S$.
We conclude that there exists a curve $\Gamma$ on $S$ that defines a 
section of $q$ and a divisor $W$ that is vertical for $q$ such that
$$
K_S+\delta \,=\,W\,+\,\Gamma\,+\,q^\ast D_a,\mbox{ \qquad }\pi^\ast\Gamma\,=\,2\widehat{C}.
$$
Applying Riemann--Roch to  $E:=\Gamma+q^\ast D_a$ on $S$ 
we compute
$$
a+1\,\geq\,h^0(S,\OO_S(E))\,\geq\,\chi(\OO_S(E))\,=\,
\chi(\OO_S)\,+\,\frac{1}{2}(E^2-EK_S).
$$
Denote by $b$ the genus of $B$.
We have $\Gamma^2+\Gamma K_S = 2b-2$, hence
$E^2-EK_S+4a=2\Gamma^2-2(b-1)+4a$.
From the previous inequality we get
$a+1\geq 1-b+\Gamma^2+(1-b)+2a$ and thus
$\Gamma^2\leq 2b-1-a$.
Coming back to $X$ we deduce
$$
 C^2\,=\,\widehat{C}^2\,\leq\, b-\frac{1}{2}(a+1)\,,
$$
whence $CK_X\geq b-2+\frac{1}{2}(a+1)$
by the adjunction formula.
Finally, we obtain
$$
K_X^2\,\geq\,2C K_X\,+\,2a\,\geq3a\,-\,3\,\geq\,3p_g\,-\,6
$$
and are done.
\qed\medskip

\begin{Remark} 
 As in \cite[Remarque 5.4]{bea} it follows that if the canonical
 map of a minimal surface of general type with $\chi(\OO_X)\geq2$ 
 is composed with a pencil over a curve of genus $b\geq1$ 
 then
 $$
   K_X^2\,\geq\,3p_g\,+\,5(b-1)
 $$
 holds true.
\end{Remark}

\subsection*{Second case: the canonical image is two-dimensional}
Let us first slightly extend the notion of hyperelliptic curves.
We refer to \cite[Section 1]{lie2} for details.

\begin{Definition}
 \label{def:hyperelliptic}
 A reduced and irreducible curve $C$, proper over an algebraically closed field,
 is called {\em hyperelliptic} if $p_a(C)\geq2$ and if there exists a morphism
 of degree $2$ from $C$ onto $\PP^1$.
\end{Definition}

Following \cite[Th\'eor\`eme 5.5]{bea} and \cite[Lemma 2]{horquint},
we now extend several results,
which are due to Noether, Castelnuovo, Horikawa, and Beauville
over the complex numbers, 
to arbitrary characteristic.

\begin{Theorem}
 \label{twodimensionalimage}
 Let $X$ be a minimal surface of general type and assume that the image
 of the canonical map is a surface.
 Let $\pi:\widetilde{X}\to X$ be a resolution of indeterminacy of the
 canonical linear system and write
 $$
    |\pi^\ast K_X| \,=\, |L| \,+\, F
 $$
 where $F$ denotes the fixed part.
 Then {\em Noether's inequality}
 $$
   K_X^2\,\geq\,L^2\,\geq\,2p_g\,-\,4
 $$
 holds true.
 Moreover,
 \begin{enumerate}
  \item If $L^2=2p_g-4$ then a general member of $|L|$ is a reduced and irreducible
    Gorenstein curve, which is hyperelliptic.
  \item If $K_X^2=L^2$ then $|K_X|$ is basepoint free.
  \item If $|L|$ defines a birational morphism then {\em Castelnuovo's inequality}
  $$
   K_X^2\,\geq\,L^2\,\geq\,3p_g\,-\,7
  $$ 
  holds true.
 \end{enumerate}
\end{Theorem}

\prf
By Bertini's theorem \cite[Th\'eor\`eme I.6.10]{jou},
a general member $C$ of $|L|$ is reduced and irreducible.
Being a divisor on a smooth surface, $C$ is Gorenstein.

Then we can conclude as in the proof of \cite[Lemma 2]{horquint}
that $K^2\geq L^2\geq2p_g-4$ and that $C$ is hyperelliptic 
(in the sense of Definition \ref{def:hyperelliptic}) if
$L^2=2p_g-4$ holds, cf. also the proof of
\cite[Theorem 2.3]{lie2}.

In order to prove assertion (2), we note that
\begin{equation}
 \label{h2vanishes}
H^2(\widetilde{X},\OO_{\widetilde{X}}(K_{\widetilde{X}}+nL))\,=\,H^0(\widetilde{X},\OO_{\widetilde{X}}(-nL))\,=\,0
\mbox{\quad for all }n\geq1
\end{equation}
by Serre duality.
By Riemann--Roch we have
\begin{equation}
 \label{rr}
\chi(\widetilde{X},\OO_{\widetilde{X}}(K_{\widetilde{X}}+nL))
\,-\, \chi(\widetilde{X},\OO_{\widetilde{X}})
\,=\,
\frac{n}{2}L(E+F)\,+\,\frac{n^2+n}{2}L^2,
\end{equation}
where $E$ denotes the divisor contained in ${\rm exc}(\pi)$
such that
$K_{\widetilde{X}}=\pi^\ast K_X+E$ holds true.

By \cite[Proposition I.1.7]{ek}
we have $h^1(-nK_X)=0$ for $n\geq2$ and then 
Riemann--Roch together with (\ref{h2vanishes}) yield
$$\begin{array}{lcl}
  {\displaystyle \frac{n(n+1)}{2} K_X^2}  &=& 
  h^0(X,\OO_X( (n+1)K_X))\,-\,\chi(\OO_X)\\
  &\geq&
  h^0(\widetilde{X},\OO_X(K_{\widetilde{X}}+nL))\,-\,\chi(\OO_{\widetilde{X}})\\
  &\geq&
  \chi(\widetilde{X},\OO_X(K_{\widetilde{X}}+nL))\,-\,\chi(\OO_{\widetilde{X}}).
\end{array}$$
Combining with (\ref{rr}) we obtain
$$
K_X^2 \,\geq\, \frac{1}{n+1}L(E+F)\,+\,L^2,
$$
for $n\geq2$ and even for $n=1$ except possibly for some surfaces 
%with $\chi=1$ 
in characteristic $2$ \cite[Theorem II.1.7]{ek}.
Now, if $L^2=K^2$, we have $LE=LF=0$.
First, $LE=0$ implies that $L$ can be considered as a divisor on $X$.
Hence $\pi$ is the identity.
Second, $|K|=|L|+F$ with $LF=0$ and the connectedness of
canonical divisors implies $F=0$.
This proves that $|K|$ has no base point in case $L^2=K^2$ holds.

To prove (3), we note that the proof of \cite[Th\'eor\`eme 5.5]{bea}
works in arbitrary characteristic.
\qed\medskip

\subsection*{Castelnuovo's inequality and Beauville's results}
The previous results extend \cite[Th\'eor\`eme 5.5]{bea} and
\cite[Remarque 5.6]{bea} to arbitrary characteristic:

\begin{Theorem}
 \label{thm:3pg-7}
 Let $X$ be a minimal surface of general type.
 \begin{enumerate}
  \item If $K^2=3p_g-7$ then the canonical map is either a rational
    map of degree $2$ onto a ruled surface or a birational morphism.
  \item If $K^2<3p_g-7$ then the canonical map is a rational map 
    of degree $2$ onto a ruled surface.
 \end{enumerate}
\end{Theorem}

Over the complex numbers,
surfaces with $K^2=3p_g-7$ and birational canonical map
have been studied by Castelnuovo \cite{castelnuovo}.
They have been classified without the birationality assumption 
by Ashikaga, Horikawa, Konno and others, 
see \cite{ashikaga konno} for details and further references.
It seems plausible that their classification 
in positive characteristic is the same
as over the complex numbers.
We shall deal with the case $K^2=3p_g-7=5$ in 
Section \ref{sec:small} below.

\section{The geometry of $\alpha_{\cal L}$-torsors}
\label{sec:alphaltorsors}

In this section we describe a handy and useful method to construct inseparable 
covers using $\alpha_{\cal L}$-torsors.
These are inseparable versions of cyclic Galois covers.
Even if the base of the cover is smooth and ${\cal L}$ is sufficiently
ample, these covers are almost always singular.
However, the singularities are controlled by Cartier's canonical
connection.
Moreover, in the surface case we will see that in nice and generic 
situations we may assume these singularities to be Du~Val. 

In characteristic zero,
cyclic Galois covers are omnipresent for the construction and 
classification of algebraic varieties.
In positive characteristic $p$, inseparable covers are equally important and
the simplest type are morphisms that are 
locally torsors under $\mu_p$ or $\alpha_p$.
Now, the group action of a Galois morphism between normal varieties is
completely determined at their function fields.
This is in contrast to inseparable morphisms: 
a morphism that is Zariski-locally an $\alpha_p$-torsor may not give rise
to a global $\alpha_p$-action.
This leads to the notion of
$\alpha_{\cal L}$-torsors introduced in \cite[Section 0]{ek},
which we recall here for the reader's convenience:
%In this section we comment on them and discuss their
%geometry in greater detail - probably some of these results 
%are known to the experts.

Let $k$ be an algebraically closed field of positive characteristic $p$
and let $S$ be a smooth variety over $k$ of arbitrary dimension.
As usual, we denote by $F:S\to S^{(p)}$ the $k$-linear Frobenius morphism.
Since $k$ is perfect we shall identify $S$ and $S^{(p)}$ whenever
needed.

Now, given an invertible sheaf $\cal L$ on $S$ we consider
$\cal L$ and ${\cal L}^{\otimes p}$ as relative 
group schemes over $S$. 
Then Frobenius induces a short exact sequence (in the flat topology)
of group schemes over $S$
\begin{equation}
 \label{alphalsequence}
 0\,\to\,\alpha_{\cal L}\,\to\,{\cal L}\,\to\,{\cal L}^{\otimes p}\,\to\,0\,.
\end{equation}
This $\alpha_{\cal L}$ is a finite flat and infinitesimal group scheme
of relative length $p$ over $S$, cf. \cite[Section 0.(1.5)]{ek}.
Locally in the Zariski-toplogy this group scheme is isomorphic to
$\alpha_p$, whence the name.

\begin{Remark}
 In characteristic $2$, every finite inseparable morphism $X\to S$ of
 degree $2$, where $S$ is smooth and $X$ is Cohen--Macaulay is
 an $\alpha_{\cal L}$-torsor for a suitable invertible sheaf
 ${\cal L}$ by \cite[Proposition 0.1.11]{ek}.
 This is analogous to the fact that every such morphism
 of degree $2$ in characteristic $\neq2$ is Galois.
\end{Remark}

Torsors under $\alpha_{\cal L}$ are classified by cohomology classes in
$H^1_{\rm fl}(S,\alpha_{\cal L})$, where the subscript
${\rm fl}$ refers to the fact that flat topology is required to define
this cohomology group (in order to trivialize such torsors).
Taking cohomology in (\ref{alphalsequence}) we obtain
\begin{equation}
 \label{alphaples}
...\,\to\,H^0(S,\,{\cal L}^{\otimes p})\,\stackrel{\delta}{\to}\, 
H^1_{\rm fl}(S,\,\alpha_{\cal L})\,\to\,H^1(S,\,{\cal L})\,\to\,...
\end{equation}
Given an open affine set $U\subseteq S$ we have
$H^1(U,{\cal L}|_U)=0$ and thus the 
$\alpha_{\cal L}$-torsor restricted to $U$
arises as $\delta(s)$ for a section $s\in H^0(U,{\cal L}^{\otimes p})$.
Choosing a local generator $t$ of ${\cal L}$ (after possibly passing
to an open affine cover of $U$), the torsor over $U$ is given by
\begin{equation}
 \label{localdescription}
 \OO_U[t]/(t^p-f)\,\to\,\OO_U
\end{equation}
for some $f\in\OO_U$.
The elements $d(s/t^p)\in\Omega^1_U$ glue to a well-defined global
$\OO_S$-morphism
\begin{equation}
 \label{torsormorphism}
    {\cal L}^{\otimes(-p)}\,\to\,\Omega_{S/k}^1\,,
\end{equation}
whose zeroes correspond to the singularities of $S$.
The annihilator of this morphism 
%(\ref{torsormorphism}) 
is a subsheaf inside the tangent
sheaf $\Theta_{S/k}$, which is the $p$-closed foliation 
associated to the inseparable morphism $\pi$, cf. 
the discussion after \cite[Lemma I.2.5]{ek}.
Viewing this assignment of a (local) section $s$ of 
${\cal L}^{\otimes p}$ to a (local) morphism (\ref{torsormorphism})
as a function of $s$, i.e.,
varying the $\alpha_{\cal L}$-torsor we obtain the following

\begin{Proposition}
 The map
 $$
 \nabla\,:\,{\cal L}^{\otimes p}\,\to\,
 \shHom({\cal L}^{\otimes(-p)},\,\Omega_{S/k}^1)
 \,\iso\,\Omega^1_{S/k}\otimes {\cal L}^{\otimes p}
 $$
 is an integrable connection of $p$-curvature zero.
 Considering ${\cal L}$ as invertible sheaf on $X^{(p)}$,
 we may identify ${\cal L}^{\otimes p}$ on $X$ with $F^\ast{\cal L}$. 
 Then horizontal sections of $\nabla$ form an invertible sheaf
 on $X^{(p)}$ and we obtain
 $$
    (F^\ast{\cal L})^\nabla \,\iso\, {\cal L}\,.
 $$
 Thus, $\nabla$ can be identified with 
 the Cartier connection on $F^*{\cal L}$.
\end{Proposition}

\prf
For a local parameter $t\in{\cal L}$ and a section $s=f\cdot t^p$ of ${\cal L}^{\otimes p}$
with $f\in\OO_S$ we have $\nabla(s)=df\otimes t^p$.
A straight forward calculation shows that $\nabla$ is an integrable
connection of $p$-curvature zero.
Moreover, its horizontal sections % $({\cal L}^{\otimes p})^\nabla$ are
are of the form $f^p\cdot t^p$, i.e., come from ${\cal L}$, 
considered as an invertible sheaf on $X^{(p)}$.
These properties characterize $\nabla$ as Cartier's connection
\cite[Theorem 5.1]{Katz}.
\qed\medskip

Let us recall that a cyclic Galois covering $\pi:X\to S$ of order $\ell$
prime to ${\rm char}(k)$ is also given by an invertible sheaf ${\cal L}$ 
and that $\pi_*\OO_X$ as a $\OO_S$-module is isomorphic to the direct sum 
of ${\cal L}^{-i}$ where $i$ runs from $0$ to $\ell-1$.
This direct sum decomposition comes from the fact that the Galois action
decomposes $\pi_*\OO_S$ into eigensheaves, all of which turn out to
be invertible if $S$ is smooth.

On the other hand, every $\alpha_{\cal L}$-torsor defines an extension
class
$$
[\,0\,\to\,\OO_S\,\to\,{\cal E}\,\to\,{\cal L}^\vee\,\to\,0\,]\,\in\,
{\rm Ext}^1({\cal L}^\vee,\OO_S)\,,
$$
which splits after Frobenius pullback \cite[Section 0.(1.6)]{ek}.
We reinterpret (\ref{alphaples}) as 
$$
...\,\to\,H^0(S,\,F^*{\cal L})\,\stackrel{\delta}{\to}\, 
H^1_{\rm fl}(S,\,\alpha_{\cal L})\,\stackrel{e}{\to}
\,{\rm Ext}^1({\cal L}^\vee,\,\OO_S)\,\to\,
\,{\rm Ext}^1({\cal L}^{\otimes(-p)},\,\OO_S)\,\to\,
...
$$
Thus, the image of $\delta$ consists of those $\alpha_{\cal L}$-torsors,
whose associated extension class in
${\rm Ext}^1({\cal L}^\vee,\OO_S)$ is trivial, i.e., torsors of the form
$$
X\,\iso\,\uSpec\,\bigoplus_{i=0}^{p-1}{\cal L}^{\otimes(-i)}\,\to\,S\,
$$
and where the multiplication map on the $\OO_S$-algebra $\OO_X$
is {\em globally} given by a section
$s\in H^0(S,{\cal L}^{\otimes p})$.
Here is an important example where the extension class of an
$\alpha_{\cal L}$-torsor is non-trivial:

\begin{Example}
 Let $E/k$ be an elliptic curve.
 Then $F:E\to E^{(p)}$ is an $\alpha_{\cal L}$-torsor
 for ${\cal L}=\OO_{E^{(p)}}$ corresponding to a
 class 
 $\eta\in H^1_{\rm fl}(E^{(p)},\alpha_{{\cal L}})$.
 It satisfies 
 $$
   0\, \neq\, e(\eta)\,\in\, {\rm Ext}^1(\OO_{E^{(p)}},\,\OO_{E^{(p)}}) 
    \,\iso\, H^1(E^{(p)},\,\OO_{E^{(p)}}),
 $$
 for otherwise the torsor would be trivial.
 Depending on whether $E$ is ordinary or supersingular,
 $F$ can also be interpreted as $\mu_p$- or $\alpha_p$-torsor,
 respectively.
\end{Example}

Let us now study $\alpha_{\cal L}$-torsors of the
form $\delta H^0(S,{\cal L}^{\otimes p})$, still assuming
$S$ to be smooth.
Thus, let $s$ be a global section of ${\cal L}^{\otimes p}$. 
Then $\delta(s)$ gives rise to
an $\alpha_{\cal L}$-torsor $\pi:X\to S$.
We have seen above that the singularities of $X$
correspond to the zeroes of $\nabla(s)$ where $\nabla$
is the Cartier connection on ${\cal L}^{\otimes p}$.

From the local description (\ref{localdescription}) it follows that
$X$ is Cohen--Macaulay since $S$ is smooth.
As explained in \cite[Section 0.(1.6)]{ek}, every non-trivial
$\alpha_{\cal L}$-torsor over a smooth base is automatically
reduced.
Clearly, a torsor of the form $\delta(s)$ is non-trivial if and
only if $s$ is not equal to $w^p$ for some
$w\in H^0(S,{\cal L})$.

If $\nabla(s)$ does not vanish along a divisor on $S$ then
the associated $\alpha_{\cal L}$-torsor $\delta(s)$ is
regular in codimension one and thus normal by Serre's criterion.

Note that $\pi:X\to S$ is everywhere ramified.
Thus, despite the description (\ref{localdescription}), the 
zero divisor of $s$ does not have a geometric interpretation as
branch divisor.
On the other hand, let $D\subset S$ be a prime divisor contained
in the zero divisor of $s$.
If $X$ is normal and $E$ denotes the reduced inverse image of $D$ then
$\pi_\ast E=D$ and $\pi^*D=pE$.
%$$
%   \pi_\ast E\,=\, D\mbox{ \quad and \quad }\pi^*D \,=\, pE\,.
%$$
%i.e., $D$ is not an integral subvariety for the $p$-closed vector
%field associated to $\pi$, cf. \cite[Proposition 1]{rs}.
\medskip

Let us finally assume that $S$ is a smooth projective surface
and $X\to S$ is an $\alpha_{\cal L}$-torsor of the form $\delta(s)$
for some $s\in H^0(S,{\cal L}^{\otimes p})$.
If $X$ is normal then $\nabla(s)$
has only isolated zeroes.
More precisely, $\nabla(s)$ will have 
$c_2(\Omega_S^1\otimes {\cal L}^{\otimes p})$
zeroes, counted with multiplicity.
Thus, in general $X$ will not be smooth, even if ${\cal L}^{\otimes p}$
is sufficiently ample.
On the other hand, the occurring singularities are not too
bad as the following Bertini-type result shows.

\begin{Theorem}
 \label{bertini}
 Let $S$ be a smooth projective surface over an algebraically closed
 field of characteristic $p>0$.
 Let ${\cal L}$ be an invertible sheaf on $S$ such that
 \begin{enumerate}
  \item $\forall x\in S$ there exist two sections 
   $s_1,s_2\in H^0(S,{\cal L}^{\otimes p})$ meeting transversally in $x$
  \item $\forall x\in S$ there exists a section 
   $s\in H^0(S,{\cal L}^{\otimes p})$ with strict normal crossings in $x$.
 \end{enumerate}
 Then there exists an open and dense subset $U\subseteq H^0(S,{\cal L}^{\otimes p})$
 such that for all $s\in U$ the associated $\alpha_{\cal L}$-torsor
 $\delta(s)$ is normal with precisely $c_2(\Omega_S^1\otimes {\cal L}^{\otimes p})$
 Du~Val singularities of type $A_{p-1}$ and no further singularities.
\end{Theorem}

\prf
The first part is similar to the proof \cite[Theorem II.8.18]{hart}
of Bertini's theorem:
For a closed point $P\in S$ we choose an open affine neighborhood
$P\in U$ and (after possibly shrinking $U$) 
a section $f_0\in H^0(U,{\cal L})$ with $f_0(P)\neq0$.
Then we define
$$\begin{array}{ccccc}
\varphi_P &:& H^0(S,\,{\cal L}^{\otimes p})&\to&\OO_{S,P}/\idealm_P^3\\
 & & s &\mapsto& \frac{s}{f_0^p}
\end{array}
$$
viewing the right hand side as $6$-dimensional 
$k$-vector space.
Using the natural embedding of $k$ into $\OO_{S,P}$ we form
the quotient map $\overline{\varphi}_P=\varphi_P\mod k$
to the $5$-dimensional vector space
$\OO_{S,P}/(k\oplus\idealm_P^3)$.

After choosing generators $x,y\in\idealm_P$ we may form the
Taylor expansion with respect to $x,y$,
which, by our particular choice of local parameter $f_0^p$
of ${\cal L}^{\otimes p}$, yields Cartier's connection
$$
\begin{array}{ccl}
\nabla(\varphi_P(s)) &=& 
\left( 
\frac{\partial\widetilde{s}}{\partial x}(P)\cdot x\,+\,
\frac{\partial^2\widetilde{s}}{\partial x^2}(P)\cdot x^2\,+\,
\frac{\partial^2\widetilde{s}}{\partial x\partial y}(P)\cdot xy\,+\,...
\right)\cdot dx\\
 &+&
\left(
\frac{\partial\widetilde{s}}{\partial y}(P)\cdot y\,+\,
\frac{\partial^2\widetilde{s}}{\partial y^2}(P)\cdot y^2\,+\,
\frac{\partial^2\widetilde{s}}{\partial x\partial y}(P)\cdot xy\,+\,...
\right)\cdot dy\,.
\end{array}
$$
Using the assumptions on the existence of sections of
$H^0(S,{\cal L}^{\otimes p})$ it is easy to see that the
image of $\overline{\varphi}_P$ is at least
$3$-dimensional.
Thus, $\ker({\overline{\varphi}_P})$ is of codimension at
least $3$ in $H^0(S,{\cal L}^{\otimes p})$.

On the other hand, a section $s\in\ker(\overline{\varphi}_P)$
has the property that $\nabla(s)$ has a zero of multiplicity
at least two in $P$.
Moreover, if the zeroes of $\nabla(s)$ are not isolated around $P$
then also $s\in\ker(\overline{\varphi}_P)$ holds true.

Considering
$$
 B \,:=\, \{ (P,s)\,|\,s\in\ker(\overline{\varphi}_P)\}\,\subseteq\,
 S\times H^0(S,\,{\cal L}^{\otimes p})
$$
with its two projections ${\rm pr}_1$, ${\rm pr}_2$ to $S$ and
$H^0(S,{\cal L}^{\otimes p})$, a dimension count
shows that ${\rm pr}_2(B)$ is a proper subset of 
$H^0(S,{\cal L}^{\otimes p})$.
Since we assumed $k$ to be algebraically closed,
there exist sections $s\in H^0(S,{\cal L}^{\otimes p})$
where the zeroes of $\nabla(s)$ are isolated of multiplicity one.
Since both properties are open it follows in fact that 
there exists an open and dense subset 
$U\subseteq H^0(S,{\cal L}^{\otimes p})$
of such sections.

Finally, locally analytically around a zero of $\nabla(s)$ 
the torsor is given by $z^p-xy$, which defines a Du~Val singularity 
of type $A_{p-1}$.
\qed\medskip

\begin{Remarks}
 We note that
 \begin{enumerate}
  \item The conditions on ${\cal L}^{\otimes p}$ are fulfilled whenever
     ${\cal L}$ is very ample.
  \item If $V\subseteq H^0(S,{\cal L}^{\otimes p})$ fulfills the assumptions
     then the conclusion also holds for a generic section of $V$.
  \item For more on generic  $\alpha_{\cal L}$-torsors, we refer to
     \cite[Chapter VI]{ek gauge}. Since I was unaware of this reference
     until pointed out to me by Torsten~Ekedahl, there are slight
     overlaps with my exposition.
 \end{enumerate}
\end{Remarks}

In \cite[I.(1.10)]{ek} the notion of $\alpha_{\cal L}$-torsor is
generalized to $\alpha_s$-torsors, which coincides with
the notion of Artin-Schreier covering of simple type in
\cite{tak}.
For a global section $s\in H^0(S,{\cal L}^{\otimes(p-1)})$ 
we define $\alpha_s$ to be the relative group scheme over
$S$, which is the kernel of 
$$
0\,\to\,\alpha_s\,\to\,{\cal L}\,\stackrel{F-s}{\to}\,{\cal L}^{\otimes p}\,\to\,0.
$$
For $s=0$ we obtain $\alpha_{\cal L}$ but if $s\neq0$ we obtain 
a group scheme, which is \'etale over $S$ outside the zero divisor
of $s$.
As before, these torsors are classified by cohomology 
classes $H^1_{\rm fl}(S,\alpha_s)$ and to every $\alpha_s$-torsor
one obtains an extension class in ${\rm Ext}^1({\cal L}^\vee,\OO_S)$.
Whenever this extension class is trivial the torsor is
globally given as
$$
\uSpec \OO_S[z]/(z^p-sz+t) \,\to\, S
$$
for some global section $t\in H^0(S,{\cal L}^{\otimes p})$.
Note that this morphism is branched over the zero divisor of $s$.
However, note that if this divisor is smooth this does not necessarily
imply that the torsor itself is smooth.
On the other hand, only singularities of type $A_{p-1}$ can occur
if the zero divisor of $s$ is smooth \cite[Theorem 2.1]{tak}.
Let us strengthen this result slightly, whose proof is analogously
to that of Theorem \ref{bertini}, which is why we leave the proof
to the reader.

\begin{Theorem}
 \label{bertini artin-schreier}
 Let $S$ be a smooth projective surface over an algebraically closed
 field of characteristic $p>0$.
 Let ${\cal L}$ be an invertible sheaf on $S$ such that
 \begin{enumerate}
  \item the sheaf ${\cal L}$ is generated by global sections,
  \item $\forall x\in S$ there exist two sections 
   $s_1,s_2\in H^0(S,{\cal L}^{\otimes(p-1)})$ meeting transversally in $x$
 \end{enumerate}
 Then there exists an open and dense subset 
 $U\subseteq H^0(S,{\cal L}^{\otimes(p-1)}) \oplus H^0(S,{\cal L}^{\otimes p})$
 such that for all $(s,t)\in U$ the associated split $\alpha_s$-torsor
 is smooth.
\end{Theorem}

\begin{Remark}
 The conditions are fulfilled whenever ${\cal L}$ is very ample.
\end{Remark}

\section{Surfaces on the Noether Lines}
\label{sec:noether}

From now on we study surfaces on the Noether lines.
In this section we prove their basic structure result except
for some cases with small $p_g$, with which we deal 
in Section \ref{sec:small}.
The surfaces of this section are double covers
of rational surfaces via their canonical map.
Next, we establish existence.
Moreover, in characteristic $2$ we prove that all classes
exist with separable as well as inseparable canonical maps.
Finally, we show that they are algebraically simply 
connected and have reduced Picard schemes.

Let us first recall from
\cite[Theorem 2.1]{lie2} 
that minimal surfaces of general type fulfill 
{\em Noether's inequality}
$$
K_X^2\,\geq\,2p_g\,-\,4\,.
$$
Over the complex numbers,
surfaces with $K^2=2p_g-4$ have been
described by Enriques in \cite[Capitolo VIII.11]{enr} and a detailed
analysis of surfaces with $K^2\leq 2p_g-3$
has been carried out by Horikawa in
\cite{horikawa} and \cite{horikawa2}.
In positive characteristic, the description
of surfaces with $K^2=2p_g-4$ is the same as in characterstic
zero \cite{lie2}.

\begin{Definition}
  A minimal surface of general type is called
  \begin{enumerate}
   \item an {\em even Horikawa surface} if $K^2=2p_g-4$, and
   \item an {\em odd Horikawa surface} if $K^2=2p_g-3$.
  \end{enumerate}
  A surface belonging to one of these classes is said to lie on the
  {\em Noether lines}.
\end{Definition}

We start with a preliminary analysis of their canonical maps, which is the
key to their classification.

\begin{Proposition}
 \label{canonicalmap}
 Let $X$ be a surface on the Noether lines.
 \begin{enumerate}
  \item 
   If $X$ is an even Horikawa surface then $p_g\geq3$ and 
   the canonical linear system is basepoint free.
  \item
   If $X$ is an odd Horikawa surface with $p_g\geq5$ then the canonical  
   linear system has a unique base point, whose indeterminacy can be resolved 
   after one blow-up.
 \end{enumerate}
 In all these cases the canonical map is generically finite of degree $2$
 onto a surface of degree $p_g-2$ in $\PP^{p_g-1}$.
\end{Proposition}

\prf
If $X$ is an even Horikawa surface then the canonical map is basepoint free
and a general member of the canonical system is a (possibly singular) 
hyperelliptic curve by Theorem \ref{twodimensionalimage}.
The canonical image is a surface which follows from 
Theorem \ref{pencil} and $p_g\geq3$.
Also, the canonical map cannot be birational since the general canonical 
divisor is hyperelliptic and so $\omega_X$ does not restrict
to a very ample linear system on it.
%, cf. also the proof of \cite[Theorem 2.3]{lie2}.
Furthermore we know %that
$$
   K_X^2\,=\,L^2\,\geq\,\deg\phi_1\cdot\,\deg \phi_1(X)\,.
$$
But $\phi_1(X)$ is a surface inside $\PP^{p_g-1}$ and thus has degree at 
least $p_g-2$.
This leaves only the possibility $\deg\phi_1=2$ and $\deg\phi_1(X)=p_g-2$.

Now, let $X$ be an odd Horikawa surface with $p_g\geq5$.
Again, by Theorem \ref{pencil} and Theorem \ref{twodimensionalimage} the image of
the canonical map is a surface and the canonical map cannot be birational.
The inequality
$$
  2p_g-3 \,=\, K^2_X \,\geq\, L^2 \,=\, \deg\phi_1\cdot \deg\phi_1(X) \,\geq\,2\cdot(p_g-2)
$$
has the only solution $K_X^2=L^2-1$, $\deg\phi_1=2$ and $\deg\phi_1(X)=p_g-2$.
In particular, the canonical system has a unique basepoint.
\qed\medskip

\begin{Remark}
 \label{delpezzo}
 A reduced and irreducible surface in $\PP^n$ that spans its ambient space has
 degree at least $n-1$. 
 Surfaces of this minimal degree have been classified by del~Pezzo.
 For $n\geq3$, the homogeneous ideal defining such a surface of degree $n-1$ in $\PP^n$ 
 is generated by $\frac{1}{2}(n-1)(n-2)$
 quadratic polynomials.
 These polynomials correspond 
 to the $2\times2$ minors of a certain $(n-1)\times2$ matrix 
 for Hirzebruch surfaces and cones over rational normal curves, 
 and the $2\times2$ minors of a certain symmetric $3\times3$ matrix 
 for the Veronese surface in $\PP^5$, namely
 $$
 \left( 
 \begin{array}{lcllcl}
  x_{0,0}&...&x_{0,a_0-1}&x_{1,0}&...&x_{1,a_1-1}\\
  x_{0,1}&...&x_{0,a_0}&x_{1,1}&...&x_{1,a_1}
 \end{array}
 \right)
 \mbox{ \quad and \quad }
 \left(
 \begin{array}{ccc}
  x_0&x_1&x_2\\
  x_1&x_3&x_4\\
  x_2&x_4&x_5
 \end{array}
 \right)\,.
 $$
 We refer to \cite{eh} for details and further references.
\end{Remark}

Now, it is not so difficult to show that the 
first Betti number of a surface on the Noether lines is zero.
Although this implies that the Picard scheme is zero-dimensional
we could still have the possibility of 
a non-trivial $H^1(X,\OO_X)$ caused by a non-reduced Picard scheme.
We start with an argument from the proof of \cite[Theorem 14]{bom}.

\begin{Lemma}
 \label{bombieri}
  Let $X$ be a surface on the Noether lines with $h^1(\OO_X)=0$.
  Then $X$ is algebraically simply connected, i.e., its \'etale
  fundamental group is trivial.
\end{Lemma}

\prf
Let $X$ be a surface with $h^1(\OO_X)=0$ and
$\widehat{X}\to X$ be an \'etale cover of degree $m$.
Then we compute $\chi(\OO_{\widehat{X}})=m\chi(\OO_X)$ and
$K_{\widehat{X}}^2=m K_X^2$.
Using Noether's inequality we obtain
$$
m(1\,+\,p_g(X))\,=\,1\,-\,h^1(\OO_{\widehat{X}})\,+\,p_g(\widehat{X})\,\leq\,1\,+\,p_g(\widehat{X})
\,\leq\,\frac{1}{2}K_{\widehat{X}}^2\,+\,3=\frac{m}{2}K_X^2\,+\,3\,.
$$
If $X$ is a surface on the Noether lines then this inequality holds for
$m=1$ only.
Thus, every \'etale cover of $X$ is trivial and hence $X$ is algebraically
simply connected.
\qed\medskip

\begin{Proposition}
  \label{h1vanishes}
  Let $X$ be a surface on the Noether lines,
  which is 
  \begin{enumerate}
   \item an even Horikawa surface or
   \item an odd Horikawa surface with $p_g\geq5$.
  \end{enumerate}
  Then $b_1(X)=h^1(\OO_X)=0$ and the Picard scheme of $X$ is
  reduced.
  Moreover, $X$ is algebraically simply connected. 
\end{Proposition}

\prf
An even Horikawa surface fulfills $p_g\geq3$, i.e., we 
may assume $p_g\geq3$ in any case.
Then the image of the canonical map is a surface by 
Theorem \ref{pencil} in case $p_g\geq4$ or 
in general by Noether's inequality \cite[Theorem 2.1]{lie2}.

From Proposition \ref{canonicalmap} we know 
that the canonical image is a surface
of degree $p_g-2$ inside $\PP^{p_g-1}$.
By Remark \ref{delpezzo} we can find a basis of
the $p_g$-dimensional vector space 
$H^0(X,\omega_X)$ such that the quadratic relations
among them are precisely the ones given
by $2\times2$ minors described there.
Using the multiplication map $|K_X|\times|K_X|\to |2K_X|$
and the quadratic relations we get
$$
  p_2(X) \,:=\, h^0(X,\,\omega_X^{\otimes2}) \,\geq\, 3p_g-3\,.
$$
Together with Riemann-Roch we find
\begin{equation}
\label{mumford}
K_X^2+\chi(\OO_X)+h^1(X,\,\omega_X^\vee) \,=\,
p_2\,\geq\,
3p_g-3\,.
\end{equation}
We have $h^1(\omega_X^\vee)\leq1$ by \cite[Theorem II.1.7]{ek}
and this dimension is zero except for a certain
class of surfaces with $\chi(\OO_X)\leq1$ in characteristic $2$.
However, (\ref{mumford}) gives us already $h^1(\OO_X)\leq2$ for our
surfaces.
In particular, they fulfill $\chi(\OO_X)\geq2$ and hence 
$h^1(\omega_X^\vee)=0$.
Applying (\ref{mumford}) once more we obtain $h^1(\OO_X)=0$ for even
Horikawa surfaces and $h^1(\OO_X)\leq1$ for odd Horikawa surfaces
with $p_g\geq5$.

We have to exclude the possibility of odd Horikawa surfaces with
$p_g\geq5$ and $h^1(\OO_X)=1$.
In this case we would have $p_2=3p_g-3$, implying that the bicanonical
map of $X$ is the composition of the canonical map followed by
the second Veronese map.
By Proposition \ref{canonicalmap}, the canonical map of such a surface 
has a basepoint and hence so has the bicanonical map.
However, this contradicts \cite[Theorem 26]{sb}.

Hence $h^1(\OO_X)=0$ holds true in all cases, and the surfaces are 
algebraically simply connected by Lemma \ref{bombieri}.
\qed

\begin{Remarks}
 Let us note that 
 \begin{enumerate}
  \item 
   Over the complex numbers, 
   these surfaces are even topologically simply
   connected by
   \cite[Theorem 3.4]{horikawa} and \cite[Theorem 4.8]{horikawa2}.
  \item 
   There do exist minimal surfaces of general type with 
   $K^2=h^{01}=p_g=1$ and $b_1=0$ in characteristic $5$,
   namely {\em non-classical Godeaux surfaces}, see
   \cite{mir} and \cite{lie godeaux}.
   Thus, not ''too far away`` from the Noether lines non-reduced
   Picard schemes show up.
%  \item 
%   This proposition supersedes and simplifies 
%   the explicit classifciation of even Horikawa surfaces obtained in  
%   \cite[Proposition 3.7]{lie2}.
 \end{enumerate}
\end{Remarks}

From del~Pezzo's analysis of surfaces of minimal degree 
and Proposition \ref{canonicalmap} we get an explicit 
description of our surfaces.

%For even Horikawa surfaces we refer to \cite[Theorem 3.3]{lie2},
%but we decided to include the statements below
%for the reader's convenience.
%Thus, we only discuss odd Horikawa surfaces 
%with $p_g\geq5$ in detail in this section.
By Proposition \ref{canonicalmap}, 
double covers play a central role.
We refer to Appendix \ref{sec:appendix}
for a characteristic-free discussion.
In particular, the notion of an invertible sheaf associated with 
a double cover is introduced there.
We remind the reader that the notion of branch divisor is not well-defined
for inseparable morphisms.
% in characteristic $2$.

For an integer $d\geq0$,
we let $\FF_d$ be the Hirzebruch surface 
$\PP_{\PP^1}(\OO_{\PP^1}\oplus\OO_{\PP^1}(d))$.
This $\PP^1$-bundle over $\PP^1$ has a section $\Delta_0$ with
self-intersection $\Delta_0^2=-d$, which is unique if
$d$ is positive.
We denote by $\Gamma$ the class of a fiber of this $\PP^1$-bundle.

\begin{Theorem}
 \label{thm:evenstructure}
  Let $X$ be an even Horikawa surface.
  Then the canonical linear system is a morphism and its canonical image
  $S=\phi_1(X)$ is a surface of degree $p_g-2$ in $\PP^{p_g-1}$. 
  Then there are the following possibilities:
  \begin{enumerate}
   \item $S$ is smooth and we have a factorization
   $$\begin{array}{ccccccccc}
 × \phi_1&:&X &\to& \Xcan &\stackrel{\pi}{\to}&S&\into&\PP^{p_g-1}\,.
   \end{array}$$
   The map $\pi$ is a finite flat morphism of degree $2$ and
   $\Xcan$ has at worst Du~Val singularities.
   If $\cal L$ denotes the invertible sheaf associated with $\pi$, then
   \begin{itemize}
    \item[-] $S\iso\PP^2$, $p_g=3$ and ${\cal L}\iso\OO_{\PP^2}(4)$
    \item[-] $S\iso\PP^2$, $p_g=6$ and ${\cal L}\iso\OO_{\PP^2}(5)$
    \item[-] $S\iso\FF_d$, $0\leq d\leq\min\{p_g-4, \frac{1}{2}(p_g+2)\}$, $p_g-d$ is even and 
      ${\cal L}\iso\OO_{\FF_d}(3\Delta_0+\frac{1}{2}(p_g+2+3d)\Gamma)$.
   \end{itemize}
   \item $S$ is the cone over a rational normal curve of degree $p_g-2$
     in $\PP^{p_g-1}$ with $4\leq p_g\leq 6$.
     There exists a partial desingularization $X'$ of $\Xcan$ such that 
     $\phi_1$ factors as
     $$\begin{array}{ccccccccccc}
 × \phi_1&:&X &\to& X' &\stackrel{\pi}{\to}&\FF_{p_g-2}\to&S &\into&\PP^{p_g-1}\,,
     \end{array}$$
     where $\pi$ is finite and flat of degree $2$. 
     If $\cal L$ denotes the invertible sheaf associated with $\pi$, then
     ${\cal L}\iso\OO_{\FF_{p_g-2}}(3\Delta_0+(2p_g-2)\Gamma)$.
  \end{enumerate}
\end{Theorem}

\prf
See \cite[Theorem 3.3]{lie2} and \cite[Proposition 4.2]{lie2}.
\qed

\begin{Theorem}
 \label{thm:oddstructure}
  Let $X$ be an odd Horikawa surface with $p_g\geq5$.
  Then the canonical linear system has a unique base point 
  and its canonical image
  $S=\phi_1(X)$ is a surface of degree $p_g-2$ in $\PP^{p_g-1}$. 
  We denote by $\widetilde{X}$ the blow-up of $X$ in this base-point and
  obtain the following possibilities for the induced morphism:
  \begin{enumerate}
   \item Case $A$: $S$ is smooth and we have a factorization
   $$\begin{array}{ccccccccc}
 × \widetilde{\phi}_1&:&\widetilde{X} &\to& \widetilde{X}' &\to&S&\into&\PP^{p_g-1}\\
   &&&\searrow&\uparrow&&\uparrow\scriptstyle\nu\\
   &&&&\widetilde{X}''&\stackrel{\pi}{\to}&\widetilde{S} 
   \end{array}$$
   where the upper row is the Stein factorization of $\widetilde{\phi}_1$.
   The morphism $\nu$ is the blow-up of $S$ in two (possibly infinitely near) 
   points $x,y$ lying on a fiber $\Gamma$ and 
   with exceptional divisors $E_x$, $E_y$.
   The map $\pi$ is a finite flat morphism of degree $2$ and
   $\widetilde{X}''$ has at worst Du~Val singularities.
   If $\cal L$ denotes the invertible sheaf associated with $\pi$, then
   \begin{center}
      $\widetilde{S}\iso\widetilde{\FF}_d$ \quad and\quad 
      ${\cal L}\iso\OO_{\widetilde{\FF}_d}(3\nu^\ast\Delta_0+\frac{1}{2}(p_g+4+3d)\nu^\ast\Gamma-2E_x-2E_y)$
   \end{center}
   where  $0\leq d\leq\min\{p_g-4, \frac{1}{2}(p_g+1)\}$ and $p_g-d$ is even.
%,  as well as $(p_g,d)\neq(5,3)$.
  \item Case $A'$: $S$ is the cone over a rational normal curve.
    Then $p_g=5$ and $S\iso\FF_3$, but $S\to\PP^4$ is no longer an embedding but
    contracts $\Delta_0$. 
    In this case, we have
   \begin{center}
      $\widetilde{S}\iso\widetilde{\FF}_3$ \quad and\quad 
      ${\cal L}\iso\OO_{\widetilde{\FF}_3}(3\nu^\ast\Delta_0+9\nu^\ast\Gamma-2E_x-2E_y)$.%,
   \end{center}
%    which corresponds to $(p_g,d)=(5,3)$ in case $A$.
  \item Case $B_1$: same as in case $A$ but only for $S\iso\FF_1$ and $p_g=5$.
    Then $\nu$ is the blow-up of $S$ in a point $x\in\Delta_0$ and
   \begin{center}
      $\widetilde{S}\iso\widetilde{\FF}_1$ \quad and\quad 
      ${\cal L}\iso\OO_{\widetilde{\FF}_1}(4\nu^\ast\Delta_0+5\nu^\ast\Gamma-2E_x)$.
   \end{center}
  \item Case $B_2$: same as in case $A$ but only for $S\iso\FF_2$ and $p_g=6$.
     Then no blow-up is needed and we have
   \begin{center}
      $S=\widetilde{S}\iso\FF_2$ \quad and\quad 
      ${\cal L}\iso\OO_{\FF_2}(4\Delta_0+7\Gamma)$.
   \end{center}
  \end{enumerate}

\end{Theorem}

\prf
For characteristic $\neq2$ this is carried out in \cite[Section 1]{horikawa2}.
In characteristic $2$ it also follows from loc.cit., only the standard
facts about double covers have to be replaced by the
characteristic $2$ description, see
Appendix \ref{sec:appendix} below.
In order to get the bounds on $p_g$ and $d$ also in characteristic $2$
we still can use \cite[Lemma 1.2]{horikawa2}, which is justified by
\cite[Lemma 4.1]{lie2}.
\qed\medskip

\begin{Remark}
 The division into cases A and B goes back to \cite[Section 1]{horikawa2}.
\end{Remark}

From the explicit description of ${\cal L}$ in 
Theorem \ref{thm:evenstructure} and Theorem \ref{thm:oddstructure} 
we get the following vanishing results, which will be important later on

\begin{Lemma}
 \label{lemma:vanishing}
 Let $\cal L$ be an invertible sheaf as in Theorem \ref{thm:evenstructure} 
 or Theorem \ref{thm:oddstructure}.
 Then 
 $$
   h^1({\cal L}^{\otimes i})\,=\,0\,\mbox{ \quad for all\quad }i\in\ZZ\,.
 $$
\end{Lemma}

\proof
For even Horikawa surfaces and odd Horikawa surfaces of type
$B_2$ this is straight forward, e.g., using
\cite[Lemma 3.5]{lie2}.
For the remaining odd Horikawa surfaces, let $\nu:\widetilde{S}\to S$ 
be the blow-up in $\{x,y\}$.
We write ${\cal L}$ as
$$
 {\cal L}\,\iso\,\nu^\ast{\cal M}\otimes\OO_{\widetilde{S}}(-2E_x-2E_y)\,.
$$
We compute
$h^1(\widetilde{S},\nu^\ast{\cal M}^{\otimes i})=h^1(S,{\cal M}^{\otimes i})=0$
for all $i\in\ZZ$ as before.
For $i\geq0$ the theorem on formal functions
gives $R^1\nu_\ast({\cal L}^{\otimes i})=0$,
and then the Grothendieck--Leray spectral sequence yields
$h^1(\widetilde{S},{\cal L}^{\otimes i})=h^1(S,\nu_\ast({\cal L}^{\otimes i}))$.
Denote by $Z_i$ the closed subscheme of $S$ 
corresponding to the ideal sheaf $\nu_\ast\OO_{\widetilde{S}}(-2iE_x-2iE_y)$. 
Taking cohomology in
$$
0\,\to\,\nu_\ast({\cal L}^{\otimes i})\,\to\,{\cal M}^{\otimes i}\,\to\,
{\cal M}^{\otimes i}|_{Z_i}\,\to\,0\,,
$$
a tedious check reveals that the boundary map 
$\delta:H^0(Z_i,{\cal M}^{\otimes i})\to H^1(S,\nu_\ast({\cal L}^{\otimes i}))$.
By the vanishing already obtained, we conclude
$h^1(\widetilde{S},{\cal L}^{\otimes i})=0$ for $i\geq0$.

Similarly, one proves  
$h^1(\widetilde{S}, {\cal L}^{\otimes -i}\otimes\omega_{\widetilde{S}})=0$
for $i\leq0$
and using Serre duality, we obtain
$h^1(\widetilde{S},{\cal L}^{\otimes i})=0$
for $i\leq0$.
\qed\medskip

Theorem \ref{thm:oddstructure} and Theorem \ref{thm:evenstructure}
give the possible list of even Horikawa surfaces and odd
Horikawa surfaces with $p_g\geq5$ in arbitrary characteristic.
It remains to prove their existence:

\begin{Theorem}
 \label{thm:existence}
 All possible cases in Theorem \ref{thm:oddstructure}
 and Theorem \ref{thm:evenstructure} do exist in arbitrary characteristic.
 Moreover, in characteristic $2$ all these types exist with inseparable 
 as well as separable canonical maps.
\end{Theorem}

\prf
In characteristic $p\neq2$, even
Horikawa surfaces are easily constructed along the lines of
\cite[Section 1]{horikawa} and odd Horikawa surfaces 
with $p_g\geq5$ along the lines of \cite[Section 1]{horikawa2}.
%Moreover, for even Horikawa surfaces in characteristic $2$ we 
%established existence in \cite[Sections 5 and 6]{lie2} using
%$p$-closed vector fields.
%Using the results of Section \ref{sec:alphaltorsors}
%we immediately reprove this result and more:

We may thus assume $p=2$.
A straight forward check shows that the sheaves ${\cal L}$
of Theorem \ref{thm:evenstructure} satisfy the assumptions of
Theorem \ref{bertini}.
This settles existence of even Horikawa surfaces
with inseparable canonical map, which we have already obtained in 
\cite[Section 5]{lie2}.

Let $\Gamma_0$ be the fiber of $\FF_d\to\PP^1$ that contains the
points $x,y$ in the statement of Theorem \ref{thm:oddstructure}. 
Consider the subspace $V\subseteq H^0(\widetilde{S},{\cal L}^{\otimes 2})$ 
consisting of those sections that vanish along the strict transform 
$\widehat{\Gamma}$ of $\Gamma_0$ on $\widetilde{S}$.
Assume furthermore that in case $A$ and $A'$ the points
$x,y$ are chosen as follows:

\begin{tabular}{lll}
 in case $A$: & $x,y\not\in\Delta_0$ & if \quad $p_g\geq3d-3$,\\
 in case $A$: & $x\in\Delta_0$, $y\not\in\Delta_0$ & if \quad $3d-3>p_g\geq2d-1$,\\
 in case $A'$: & $x\in\Delta_0$.
\end{tabular}

A straight forward, yet tedious calculation reveals
that the subspace $V$ fulfills the assumptions of 
Theorem \ref{bertini}.
Note that $\pi^*\widehat{\Gamma}=2E$ for some rational curve $E$
on the associated inseparable cover and thus $E^2=-1$.
Hence $E$ is an exceptional $(-1)$-curve as in 
characteristic $\neq2$.
This settles existence of all cases of Theorem \ref{thm:oddstructure}
with inseparable canonical map in characteristic $2$.

From Lemma \ref{lemma:vanishing} we get
${\rm ext}^1({\cal L}^\vee,\OO_X)=h^1({\cal L})=0$ and so
(\ref{associated}) splits.
Thus, the double cover $\pi$ for a Horikawa surface is of
the form
$\pi:\uSpec (\OO_{\widetilde{S}}\oplus {\cal L}^\vee)\to\widetilde{S}$.
The $\OO_{\widetilde{S}}$-structure is globally given by 
$z^2+fz+s=0$ where 
$s$ is a global section of ${\cal L}^{\otimes 2}$.
The canonical map is inseparable
if and only if $f=0$.
Since $h^0({\cal L})\neq0$ in all cases, a generic global section $f$
of ${\cal L}$ yields a Horikawa surface with separable canonical map.
We refer to \cite[Section 6]{lie2} for technical details.
\qed

\section{Moduli spaces} 
\label{sec:moduli}

In this section we give a characteristic-free description of
moduli spaces of even Horikawa surfaces with fixed canonical
image, which turn out to be irreducible and unirational.
In characteristic $2$, surfaces with 
inseparable canonical map form an irreducible, unirational, 
proper and closed subscheme inside these moduli spaces. 
We restrict ourselves to even Horikawa surfaces in order
to keep our discussion short.
The interested reader will have no difficulties working out
the corresponding statements for odd Horikawa surfaces.

We obtained the case of separable canonical maps via a deformation
argument from the inseparable ones.
This is not a coincidence as we shall see now.
To make this precise, let $S$ be a surface of minimal degree $n-2$ 
inside $\PP^{n-1}$. 
Then we define the following sets
$$\begin{array}{lcl}
 \Mev(S,n) &:=& \{ \mbox{ even Horikawa surfaces with $p_g=n$ and } \\
             & & \mbox{ \quad canonical image $S$ } \}\,/\,\iso \\
 \Mevi(S,n) &:=& \{ \mbox{ surfaces in $\Mev(S,n)$ with inseparable $\phi_1$ } \} 
% \\
% \Modd(S,n) &:=& \{ \mbox{ odd Horikawa surfaces with $p_g=n$ and } \\
%             & & \mbox{ \quad canonical image $S$ } \}\,/\,\iso \\
% \Moddi(S,n) &:=& \{ \mbox{ surfaces in $\Modd(S,n)$ with inseparable $\phi_1$ } \}
\end{array}$$
Whether the set $\Mev(S,n)$ is empty or not is answered by
Theorem \ref{thm:evenstructure} and Theorem \ref{thm:existence}.
%For $n\geq5$, Theorem \ref{thm:oddstructure} and Theorem \ref{thm:existence}
%gives (non-)emptyness for $\Modd(S,n)$.
%
Here is a characteristic-free description of these spaces:

\begin{Theorem}
 \label{thm:evenmoduli}
 Whenever $\Mev(S,p_g)$ is not empty it 
 carries the structure of an irreducible and
 unirational scheme of dimension
 $$
   h^0(S,{\cal L}^{\otimes 2}) \,-\, \dim\Aut(S) \,-\, 1\,,
 $$
 where $\cal L$ is as in Theorem \ref{thm:evenstructure}.
 Whenever
 \begin{itemize}
  \item [-] $S\iso\PP^2$, which implies $p_g=3$ or $p_g=6$, or
  \item [-] $S\iso\FF_d$ and $p_g\geq 3d-2$ or $p_g=2d-2$, or
  \item [-] $S$ is not smooth and $p_g=4$ or $p_g=6$,
 \end{itemize}
 then the generic surface of $\Mev(S,p_g)$ has a smooth canonical model.

 In the remaining cases the generic canonical model has
 $A_1$-singularities.
 The number of $A_1$-singularities of the
 generic canonical model is as follows:
 $$\begin{array}{l||c|c}
   (S,\,p_g)               & (\FF_d,\,p_g),\, 3d-2\,>\,p_g\,>\,2d-2 & (\mbox{cone},\,5) \\
   \hline  
   A_1-\mbox{singularities}    & p_g\,+\,2\,-\,2d & 1 
 \end{array}$$
\end{Theorem}

\prf
Assume that $S$ is smooth and let $X\in\Mev(S,p_g)$.
From the vanishing ${\rm ext}^1({\cal L}^\vee, \OO_S)=h^1({\cal L})=0$
by Lemma \ref{lemma:vanishing} it follows that $\pi:\Xcan\to S$ is isomorphic 
to $\pi:\uSpec (\OO_S\oplus {\cal L}^\vee)\to S$.
The $\OO_S$-algebra structure is given by $z^2+fz+g=0$ where
$f\in H^0(S,{\cal L})$ and $g\in H^0(S,{\cal L}^{\otimes 2})$.

Conversely, for generic choices of global sections of ${\cal L}$ and
${\cal L}^{\otimes 2}$, the associated double cover $\pi$ will be
the canonical model of an even Horikawa surface with canonical image
$S$ inside $\PP^{p_g-1}$.
In characteristic $\neq2$ this is classical.
In characteristic $2$ the argument is as follows:
a generic global section $g\in H^0(S,{\cal L}^{\otimes 2})$ 
and $0=f\in H^0(S,{\cal L})$ give rise to an even Horikawa surface 
with the desired properties by Theorem \ref{bertini}.
But then, fixing $g$, a generic choice of $f\in H^0(S,{\cal L})$
yields again a surface with at worst Du~Val 
singularities \cite[Proposition 6.1]{lie2}, i.e., the canonical
model of an even Horikawa surface with the desired properties.
Thus, the assertion is also true in characteristic $2$.

Therefore, there exists an open and dense subset
$U\subseteq H^0(S,{\cal L})\oplus H^0(S,{\cal L}^{\otimes 2})$
mapping surjectively onto $\Mev(S,p_g)$.
Now, suppose we are given two $\OO_S$-algebras $\OO_S[z_i]/(z_i^2+f_iz_i+g_i)$
that are isomorphic over $S$.
A straight forward computation shows that such an isomorphism
is given by $z_1\mapsto a\cdot z_2+h$ for some $a\in k^\times$
and $h\in H^0(S,{\cal L})$.
Conversely, such a map defines an isomorphism over $S$.
Thus, given two isomorphic even Horikawa surfaces, there exists
an isomorphism of $S$ after which they become isomorphic as
surfaces over $S$.
Thus, $\Mev(S,p_g)$ arises from $U$ by quotienting out the indicated
actions of $\Aut(S)$, $\GG_m$ and $H^0(S,{\cal L})$.
In particular, $\Mev(S,p_g)$ is irreducible, unirational and of
the stated dimension.

In case $S$ is a cone the arguments are similar and left to the reader.

Finally, if $S=\PP^2$, $S=\FF_d$ and $p_g>3d-2$ 
or $p_g=4$ and $S$ is a cone then ${\cal L}$ and
${\cal L}^{\otimes 2}$ are very ample.
In characteristic $\neq2$ we can thus choose $0=f\in H^0(S,{\cal L})$
and $g\in H^0(S,{\cal L}^{\otimes2})$ with smooth zero divisor and
the associated double cover will be smooth.
In characteristic $2$ this conclusion still holds true by
Theorem \ref{bertini artin-schreier}.
In the remaining cases,
sections of ${\cal L}$ and ${\cal L}^{\otimes 2}$ necessarily
vanish along $\Delta_0$.
We leave it to the reader that in case $S=\FF_d$ and $p_g=3d-2$
and $S$ a cone and $p_g=6$ it is still possible to find
smooth sections whereas this is impossible in the other cases.

As to the number of $A_1$-singularities, 
let us first deal with $S=\FF_d$ and $3d-2>p_g>2d-2$.
In characteristic $p\neq2$, a generic section of 
${\cal L}^{\otimes2}$ vanishes along the union of $\Delta_0$ 
and a smooth divisor intersecting
$\Delta_0$ transversally in $p_g+2-2d$ points.
Thus, the associated double cover has the claimed number of
$A_1$-singularities.
In characteristic $2$ a similar argument and
Theorem \ref{bertini} reveal that an inseparable
cover associated to a generic section of ${\cal L}^{\otimes2}$
has at least $A_1$ singularities at $p_g+2-2d$ points
lying on $\Delta_0$ (but in fact more, see 
Theorem \ref{thm:eveninsepmoduli} below).
A generic section of ${\cal L}$ vanishes along $\Delta_0$
and a smooth divisor intersecting transversally.
This and a local computation reveals that a
generic Artin--Schreier covering still has
$p_g+2-2d$ singularities of type $A_1$.
We leave the details and the case where $S$ is a cone
to the reader.
\qed

\begin{Remark}
 Since $S$ and ${\cal L}$ are defined over the integers,
 the proof provides an explicit construction of
 $\Mev(S,p_g)$ over $\Spec\ZZ$.
 See also Theorem \ref{thm:flat} below.
\end{Remark}

We included the statement about generic canonical models being
smooth or not for the following reason:
in case the canonical map is inseparable the canonical
model has always singularities:

\begin{Theorem}
 \label{thm:eveninsepmoduli}
 Suppose that we are in characteristic $2$ and that $\Mev(S,p_g)$ is not empty. 
 Then the subset 
 $$
   \Mevi(S,p_g) \mbox{ \quad }\subseteq\mbox{ \quad }\Mev(S,p_g)
 $$ 
 forms an irreducible and unirational subscheme of codimension $h^0(S,{\cal L})$.

 Every surface of $\Mevi(S,p_g)$ has a singular canonical model, but the generic
 member has only Du~Val singularities of type $A_1$.
 The number of $A_1$-singularities of the
 generic canonical model is as follows:
 $$\begin{array}{l|cccc}
   (S,\,p_g)               & (\PP^2,\,3) & (\PP^2,\,6) & (\FF_d,\,p_g) & (\mbox{cone},\,p_g) \\
   \hline  
   A_1-\mbox{singularities}  & 43 & 73 & 12\,+\,10\,p_g & 12\,+\,10\,p_g
 \end{array}$$

 In characteristic $\neq2$, the set $\Mevi(S,p_g)$ is empty.
\end{Theorem}

\prf
The canonical map of an even Horikawa surface is generically finite
of degree $2$ and so $\Mevi(S,p_g)$ is empty in characteristic $\neq2$.
In characteristic $2$, surfaces with inseparable canonical map
correspond to those $\OO_S$-algebras $\OO_S[z]/(z^2+fz+g)$,
where $0=f\in H^0(S,{\cal L})$, keeping the notations
of the proof of Theorem \ref{thm:evenmoduli}.
By the arguments given there we find an open dense
subset of $H^0(S,{\cal L}^{\otimes2})$ mapping surjectively onto
$\Mevi(S,p_g)$.
Taking into account the actions of $H^0(S,{\cal L})$, $\GG_m$
and $\Aut(S)$ that link isomorphic surfaces we conclude as before
that $\Mevi(S,p_g)$ is irreducible, unirational and obtain the 
asserted codimension.

Since $c_2(\Omega_S\otimes {\cal L}^{\otimes2})$ is non-zero in all cases
the inseparable double covers can never be smooth, see
Section \ref{sec:alphaltorsors}.
However, by Theorem \ref{bertini} the singularities of a 
generic member of $\Mevi(S,p_g)$ are Du~Val singularities
of type $A_1$.

In order to determine the number of $A_1$ singularities
let $\pi:\Xcan\to S$ be the canonical model and its
inseparable double cover, and assume $S$ is smooth.
Let $\widetilde{\pi}:X\to\widetilde{S}$ be the canonical resolution 
of singularities as given by Proposition \ref{prop:canres}.
The singularities of $\Xcan$ all being of type $A_1$,
the number of these singularities equals the number
of blow-ups in closed points $\widetilde{S}\to S$,
which is equal to $b_2(\widetilde{S})-b_2(S)$.
On the other hand, $\widetilde{\pi}$ is finite and
purely inseparable, and $X$ and $\widetilde{S}$ are smooth,
which implies $b_2(X)=b_2(\widetilde{S})$.
Since $X$ is an even Horikawa surface we have
$b_2(X)+2=c_2(X)=12\chi(\OO_X)-K^2$ from which it
is easy to compute the stated numbers.
In case $S$ is not smooth we have to blow up once
to start with and conclude as before.
\qed\medskip

\begin{Remarks}
 We note
 \begin{enumerate}
  \item Even if every canonical
   model of $\Mev(S,p_g)$ in characteristic $2$
   is singular, the number of singularities of
   a generic member of $\Mev(S,p_g)$ is less
   than that of a generic member of $\Mevi(S,p_g)$.
  \item The explicit examples of even Horikawa surfaces with
    inseparable canonical map of 
    \cite[Sections 5 and 6]{lie2} have $D_4$-singularities
    and are thus rather special.
 \end{enumerate}
\end{Remarks}
\vfill

\pagebreak

\section{Small invariants}
\label{sec:small}

We are left with the description of odd Horikawa surfaces with $p_g\leq4$.
In this case the canonical map need no longer be generically finite of degree $2$
onto a surface of minimal degree and new cases show up:

\begin{Proposition}
  \label{small}
  Let $X$ be an odd Horikawa surface with $p_g\leq4$.
  Then we have the following cases
  \begin{enumerate}
   \item $p_g=2$, $K_X^2=1$ and the canonical map is a pencil.
   \item $p_g=3$, $K_X^2=3$, the canonical linear system has a unique base point. 
     After blowing it up we obtain a generically finite
     morphism of degree $2$ onto $\PP^2$.
   \item $p_g=3$, $K_X^2=3$, the canonical linear system has no basepoints and  
     defines a generically finite morphism of degree $3$ onto $\PP^2$.
   \item $p_g=4$, $K_X^2=5$, the canonical linear system has no basepoints and  
     defines a birational morphism onto a quintic surface in $\PP^3$.
   \item $p_g=4$, $K_X^2=5$, the canonical linear system has a unique basepoint.
     After blowing it up we obtain a generically finite
     morphism of degree $2$ onto a (possibly singular) quadric surface in $\PP^3$.
  \end{enumerate}
\end{Proposition}

\prf
If $p_g=2$ then the canonical map defines a pencil, whereas for $p_g\geq3$
the canonical image is always a surface by \cite[Theorem 2.1]{lie2}.
As in the proof of Proposition \ref{canonicalmap} the result follows
from Theorem \ref{twodimensionalimage}.
\qed

\subsection*{Surfaces with $\mathbf{p_g=2}$ and $\mathbf{K_X^2=1}$}

These surfaces have the smallest invariants possible among all surfaces on the
Noether lines.
The canonical map defines a pencil, which does
not happen for surfaces on the Noether lines with larger invariants. 
%cf. Proposition \ref{small} and Proposition \ref{canonicalmap}.
We use the bicanonical map to obtain an explicit 
classification.

\begin{Theorem}
   Let $X$ be a minimal surface of general type with $p_g=2$ and $K^2=1$. 
   Then the bicanonical map is a morphism and factors as
   $$
    \phi_{|2K_X|}\,:\,
    X\,\to\,\Xcan\,\stackrel{\pi}{\to}\,\FF_2\,\stackrel{\nu}{\to}\,\PP^3\,,
   $$ 
   where $\pi$ is a finite flat morphism of degree $2$ onto the Hirzebruch
   surface $\FF_2$ and with associated invertible sheaf 
   ${\cal L} \,\iso\,\OO_{\FF_2}(3\Delta_0\,+\,5\Gamma)$.
   The map $\nu$ contracts the $(-2)$-curve $\Delta_0$ on $\FF_2$ and maps
   $\FF_2$ onto a quadric cone in $\PP^3$.
   Moreover, $h^1(\OO_X)=0$ and the surface is algebraically
   simply connected.

   In characteristic $2$, such surfaces exist with separable as well as
   inseparable bicanonical maps.
\end{Theorem}

\prf
Ekedahl's inequality
\cite[Corollary II.1.8]{ek} yields $h^1(\OO_X)\leq1$.
Thus, we obtain $\chi(\OO_X)\geq2$ and so $h^1(\omega_X^\vee)=0$
by \cite[Corollary II.1.8]{ek}.
In particular, $h^0(\omega_X^{\otimes 2})=3$.
Arguing as in \cite[Lemma 2.1]{horikawa2}, we 
conclude that the bicanonical linear
system defines a morphism of degree $2$ onto a quadric cone in
$\PP^3$.
As in \cite[Section 2]{horikawa2} we obtain the structure result of
these surfaces.
From the vanishing $h^1({\cal L}^\vee)=0$ and
(\ref{associated}) we conclude $h^1(\OO_X)=0$
and the surface is algebraically simply connected by Lemma \ref{bombieri}.

The existence of such surfaces in characteristic $\neq2$ is shown along
the lines of \cite[Section 2]{horikawa2}.
Surfaces with inseparable canonical maps arise as
$\alpha_{\cal L}$-torsors over $\FF_2$  
and applying Theorem \ref{bertini}.
Deforming these surfaces, we obtain surfaces 
with separable bicanonical maps.
\qed

\subsection*{Surfaces with $\mathbf{p_g=3}$ and $\mathbf{K_X^2=3}$}

Here there are two cases and both do exist. 
Also, we encounter possibly inseparable canonical maps
in characteristic $3$.

\begin{Theorem}
   \label{thm:pg=k2=3}
   Let $X$ be a minimal surface of general type with $p_g=K^2=3$. 
   Then 
   \begin{enumerate}
    \item either the canonical map has a unique base point and 
      after blowing it up we obtain a morphism of degree $2$
      onto $\PP^2$.

      In characteristic $2$, such surfaces exist with separable as
      well as inseparable canonical maps.
    \item or the canonical map is a generically finite morphism of
     degree $3$ onto $\PP^2$.

     In characteristic $3$, such surfaces exist with separable as
     well as inseparable canonical maps.
   \end{enumerate}
   In any case, $h^1(\OO_X)=0$ and the surface is algebraically
   simply connected.
\end{Theorem}

\prf
By Proposition \ref{small} there can only be these two cases.

In case the canonical map has a unique base point, let 
$\widetilde{X}\to X$ be the blow-up of this base point and
the canonical map extends to a generically finite morphism
$\widetilde{\phi}_1:\widetilde{X}\to\PP^2$.
As explained in \cite[Section 2]{horikawa2}, the exceptional
$(-1)$-curve of $\widetilde{X}$ maps to a line $\ell$ in $\PP^2$
and $\widetilde{X}$ determines three points $x,y,z$ on $\ell$.
Let $q:\widetilde{\PP^2}\to\PP^2$ be the blow-up in these three 
points with exceptional divisors $E_x$, $E_y$ and $E_z$.
Then $\widetilde{\phi}_1$ factors over $\widetilde{\PP^2}$.
Let $\widetilde{X}\to X'\to\widetilde{\PP^2}$ be the Stein 
factorization.
Then $X'\to\widetilde{\PP^2}$ is a finite flat
double cover with associated invertible sheaf 
${\cal L}^\vee=q^*\OO_{\PP^2}(-5)\otimes\OO_{\widetilde{\PP^2}}(2E_x+2E_y+2E_z)$.
From the vanishing $h^1({\cal L})=0$ we
get existence as in Section \ref{sec:moduli}.
Moreover, from $h^1({\cal L}^\vee)=0$ we conclude $h^1(\OO_X)=0$
and thus, $X$ is algebraically simply connected by
Lemma \ref{bombieri}. 

In case the canonical map is a morphism of degree $3$ we
can factor it as $X\to\Xcan\to\PP^2$, where $\Xcan$ is
the canonical model and $\Xcan\to\PP^2$ is finite and flat
of degree $3$.
If $h^1(\OO_X)=0$ then $X$ is algebraically simply
connected.
Moreover, in this case one can argue as in 
\cite[Section 2]{horikawa2} to conclude that $\Xcan$ embeds
as a hypersurface given by an equation
$\psi^3+A\psi^2+B\psi+C=0$
%$$
%\psi^3\,+\,A\psi^2\,+\,B\psi\,+\,C\,=\,0
%$$
into the vector bundle  
$V(\OO_{\PP^2}(2))\to\PP^2$.
Here, $A,B,C$ are elements of ${\rm Sym}^i H^0(\omega_X)$
with $i=2,4,6$.
Existence in characteristic $\neq3$ follows as in
loc. cit.
Surfaces with inseparable canonical maps can be constructed
using $\alpha_{\cal L}$-torsors over $\PP^2$ with 
${\cal L}=\OO_{\PP^2}(2)$ and applying 
Theorem \ref{bertini}.
Deforming these surfaces, we obtain surfaces 
with separable canonical maps.

It remains to exclude the existence of surfaces
whose canonical map is a morphism of degree $3$ and 
$h^1(\OO_X)\neq0$.
As in the proof of Proposition \ref{h1vanishes} such surfaces
fulfill $h^1(\OO_X)=1$ and we obtain $\chi(\OO_X)=3$.
Then $\pi:\Xcan\to\PP^2$ is a finite flat morphism of degree $3$
and we obtain a short exact sequence
$$
0\,\to\,\OO_{\PP^2}\,\to\,\pi_\ast\OO_{\Xcan}\,\to\,{\cal E}\,\to\,0\,,
$$
where $\cal E$ is a locally free sheaf of rank $2$ on $\PP^2$.
From $K^2=3$ and $\chi(\OO_X)=3$ 
we get two equations for the Chern classes of ${\cal E}$, see
\cite[Proposition 8.1]{pardini} and \cite[Corollary 8.3]{pardini}.
A straight forward computation shows that these equations cannot
be fulfilled by any locally free sheaf of rank $2$ on $\PP^2$.
Thus $\cal E$ and $X$ do not exist.
\qed

\subsection*{Surfaces with $\mathbf{p_g=4}$ and $\mathbf{K_X^2=5}$}

Again, we have to deal with two cases.
We note that this class of Horikawa surfaces also fulfills
$K^2=3p_g-7$.
Thus, Proposition \ref{small} overlaps with the class of 
surfaces described in Theorem \ref{thm:3pg-7}.

\begin{Theorem}
   Let $X$ be a minimal surface of general type with $p_g=4$ 
   and $K_X^2=5$.
   Then 
   \begin{enumerate}
    \item either the canonical map is a birational morphism
     that embeds the canonical model $\Xcan$ as a quintic surface
     into $\PP^3$,  
    \item or the canonical map has a unique base point. 
      After blowing it up we obtain a morphism of degree $2$
      onto a possibly singular quadric surface in $\PP^3$.

      In characteristic $2$, such surfaces exist with separable as
      well as inseparable canonical maps.
   \end{enumerate}
   In any case, $h^1(\OO_X)=0$ and the surface is algebraically
   simply connected.
\end{Theorem}

\prf
By Proposition \ref{small} (or Theorem \ref{thm:3pg-7}) 
there can only be these two cases.

In the first case, Theorem \ref{thm:3pg-7} implies that $K_X$ is
very ample on $\Xcan$, embedding it as a quintic surface in $\PP^3$.
Since quintic surfaces in $\PP^3$ with at worst Du~Val singularities
fulfill $h^1=0$ and are algebraically simply connected,
the same is true for $X$.

In the second case, Proposition \ref{small} yields their structure.
Thus, we have to consider double covers of $\PP^1\times\PP^1$
with associated invertible sheaf
% (see Section \ref{sec:appendix})
${\cal L}\iso\OO_{\PP^1\times\PP^1}(3,3)$, or
of  $\FF_2$ with 
${\cal L}\iso\OO_{\FF_2}(3\Delta_0+3\Gamma)$.
From $h^1({\cal L}^\vee)=0$ and (\ref{associated}) we deduce 
$h^1(\OO_X)=0$ and hence these surfaces are algebraically 
simply connected by Lemma \ref{bombieri}.
Existence is shown as in Section \ref{sec:moduli}.
\qed

\section{Hodge degeneration and crystalline cohomology}
\label{sec:hodge}

In this section we show that 
that the Fr\"olicher spectral sequence
%from Hodge- to deRham-cohomology 
for surfaces on the Noether lines
always degenerates at $E_1$.
Moreover, we show that their Hodge- and Betti- numbers
coincide.
In particular, their crystalline cohomology groups
are torsion-free.
On the other hand, we show that there exist even Horikawa
surfaces in arbitrary large characteristics, whose
%Hodge--Witt cohomology groups are not finitely generated.
%For these surfaces the 
slope spectral sequence 
from Hodge--Witt- to crystalline cohomology 
does not degenerate at $E_1$.

\begin{Proposition}
 \label{global 1-forms}
 Let $X$ be a surface on the Noether lines.
 Then $H^0(X,\Omega_X^1)=0$.
 In particular, all global $1$-forms are $d$-closed.
\end{Proposition}

\prf
We will only deal with the cases where the canonical or bicanonical 
map is generically finite of degree $2$ and
leave quintic surfaces in $\PP^3$ 
(which satisfy $p_g=4$, $K^2=5$) and triple covers
of $\PP^2$ with $p_g=K^2=3$ to the reader.

Then, by our explicit classification there exists a birational
model $X'$ of $X$ with at worst Du~Val singularities
and a finite flat morphism $\pi:X'\to S$ of degree $2$ onto
a smooth surface $S$.
Let $\cal L$ be the associated invertible sheaf.
By Lemma \ref{lemma:vanishing}, we have 
${\rm Ext}^1({\cal L}^\vee,\OO_S)=H^1(S,{\cal L})=0$,
which implies that the
exact sequence (\ref{associated}) splits.
Thus, $X'$ embeds into the total space of the 
vector bundle 
$$\xymatrixcolsep{5pc}
\xymatrix{
 X' \ar[r]^-i \ar[rd]^\pi & V  \,=\,\uSpec \left( 
  \bigoplus_{n\geq0}{\rm Sym}^n({\cal L}^\vee) \right)  \ar[d]_q \\
 & S
}$$
%$$
%\begin{array}{ccc}
% X &\stackrel{i}{\into}& V\,=\,\uSpec \bigoplus_{n\geq0}{\rm Sym}^n({\cal L}^\vee) \\
% &\searrow&\downarrow\\
% & & S
%\end{array}
%$$
The cotangent sequence for $q$ on $V$ exhibits 
$\Omega_V^1$ as an extension of $q^\ast({\cal L}^\vee)$ 
by $q^\ast\Omega_S^1$.
In particular, $H^0(V,\Omega_V^1)=0$.
Passing to reflexive hulls, i.e., taking double
duals, the conormal sequence becomes
\begin{equation}
 \label{conormal}
 0\,\to\,{\cal N}_{X/V}^{\vee\vee}
 \,\to\,i^\ast(\Omega_V^1)\,\to\,
 (\Omega_{X'}^1)^{\vee\vee}\,\to\,0\,,
\end{equation}
where injectivity on the left follows 
from torsion-freeness.
Note that $\det((\Omega_{X'}^1)^{\vee\vee})$
and $\omega_{X'}$ are isomorphic, since they
coincide outside codimension two and are both 
reflexive sheaves on a normal variety. 
Taking determinants in (\ref{conormal}) and using 
$\omega_{X'}\iso\pi^\ast(\omega_S\otimes{\cal L})$ 
we see that ${\cal N}_{X'/V}^{\vee\vee}$
is isomorphic to $\pi^\ast({\cal L}^{-2})$.
Pushing ${\cal N}_{X'/V}^{\vee\vee}$ down to $S$,
and using 
$H^1(S,{\cal L}^{-2})=H^1(S,{\cal L}^{-3})=0$
from Lemma \ref{lemma:vanishing},
we conclude $H^1(X',{\cal N}_{X'/V}^{\vee\vee})=0$.
Thus, taking cohomology in (\ref{conormal}) we obtain 
$H^0(X',(\Omega_{X'}^1)^{\vee\vee})=0$.

Now, let $\nu:\widetilde{X}\to X'$ be a resolution of singularities
with exceptional divisor $Z$ and
and let $U$ be the smooth locus of $X'$.
By reflexivity we have 
$H^0(U,\Omega^1_U)=H^0(X',(\Omega_{X'}^1)^{\vee\vee})=0$.
On the other hand, since $\nu$ is the resolution of Du~Val 
singularities there are never global sections of
$\Omega_{\widetilde{X}}^1$ entirely supported on $Z$ (in fact,
global sections of $\Omega_{\widetilde{X}}^1$ vanish along $Z$).
%Thus, the local cohomology group $H^0_Z(\widetilde{X},\Omega_{\widetilde{X}}^1)$
%is zero and 
Thus, the restriction map
$H^0(\widetilde{X},\Omega_{\widetilde{X}}^1)\to H^0(U,\Omega_U^1)$
is injective and
we conclude $H^0(X,\Omega_X^1)=H^0(\widetilde{X},\Omega_{\widetilde{X}}^1)=0$.
\qed\medskip

Let us recall that the Fr\"olicher spectral sequence
$$
E_1^{p,q}\,=\,H^q(X,\,\Omega_X^p)\,\Rightarrow\,
\HdR{p+q}(X)
$$
links Hodge cohomology to deRham cohomology of smooth projective
varieties.
Over fields of characteristic zero, as well as for
curves and Abelian varieties over arbitrary ground
fields, this spectral sequence degenerates
at $E_1$.

Let us also recall, e.g. from \cite[II.(4.9.1)]{illusie},
the universal coefficient formula for crystalline 
cohomology:
$$
0\,\to\,\Hcris{i}(X/W)\otimes_W k\,\to\,\HdR{i}(X/k)\,\to\,
{\rm Tor}_1^W(\Hcris{i+1}(X/W),\,k)\,\to\,0\,.
$$
It shows that deRham- and Betti-numbers coincide if
and only if all crystalline cohomology groups are
torsion-free.

\begin{Theorem}
 \label{froelicher}
 Let $X$ be a surface on the Noether lines.
 Then
 \begin{enumerate}
  \item its Fr\"olicher spectral sequence degenerates at $E_1$-level, and
  \item all crystalline cohomology groups are torsion-free.
 \end{enumerate}
 In particular, the deRham- and Betti-numbers coincide.
\end{Theorem}

\prf
We have $h^{10}=h^{01}=0$ for surfaces on the Noether lines
by the classification results of the previous sections and
Proposition \ref{global 1-forms}.
Moreover, we get $h^{12}=h^{21}=0$ from Serre duality.
Already the existence of the Fr\"olicher spectral sequence
gives the inequalities $h^1_{\rm dR}\leq h^{01}+h^{10}$ and
$h^3_{\rm dR}\leq h^{21}+h^{12}$.
Thus, we obtain $h^1_{\rm dR}=h^3_{\rm dR}=0$.
Moreover,
$\sum_n (-1)h^n_{\rm dR}$ coincides with
$\sum_{i+j}(-1)^{i+j}h^{ij}$,
which implies
$h^2_{\rm dR}=h^{02}+h^{11}+h^{20}$.
These equalities prove degeneration of the 
Fr\"olicher spectral sequence at $E_1$.

From $h^1_{\rm dR}=h^3_{\rm dR}=0$ and the universal
coefficient formula for crystalline cohomology
we obtain $\Hcris{1}=\Hcris{3}=0$ and see that 
$\Hcris{2}$ is torsion-free.
\qed\medskip

In \cite{illusie} the Hodge--Witt cohomology
groups $H^j(X,W\Omega_X^i)$ are constructed.
For smooth projective varieties,
these are $W=W(k)$-modules of finite rank, 
whose torsion-groups may not be finitely generated.
The slope spectral sequence
$$
  E_1^{p,q}\,=\,H^q(X,\,W\Omega_X^p)\,\Rightarrow\,\Hcris{p+q}(X/W)
$$
links Hodge--Witt- to crystalline cohomology.
For surfaces, it degenerates at $E_1$ if and only if
$H^2(X,W\OO_X)$ is finitely generated
\cite[Corollaire II.3.14]{illusie}.
Moreover, it degenerates always modulo torsion
at $E_1$, see \cite[Th\'eor\`eme II.3.2]{illusie}.

%The proof of Theorem \ref{froelicher} shows that $\Hcris{1}(X/W)=0$
%for surfaces on the Noether lines.
%By partial degeneration of the slope spectral sequence
%\cite[Proposition II.3.11]{illusie} this implies
%$$
%  H^1(X,\,W\OO_X)\,=\,H^0(X,\,W\Omega_X^1)\,=\,0\,.
%$$
In our case, this finite generation is closely related
to the formal Brauer group:
namely, if $X$ is a surface on the Noether lines then
$h^{01}=h^{03}=0$ for dimensional reasons and
by the classification results of the
previous sections.
Thus, the functor on Artin-algebras over $k$
$$\begin{array}{ccccc}
   {\rm Br} &:& A &\mapsto& {\rm ker}\left(\, 
 H^2(X\times_k A,\,\GG_m)\,\to\,H^2(X,\,\GG_m)
   \,\right)
  \end{array}
$$
is pro-representable by a smooth formal group $\fBr(X)$
of dimension $p_g=h^{02}$, the formal Brauer group.
Moreover, its Dieudonn\'e module of typical curves
can be identified with $H^2(X,W\OO_X)$.
In particular, it is finitely generated
if and only if $\fBr(X)$ is of finite height, i.e.,
has no unipotent part \cite[Section II.4]{artin mazur}.

In order to show that $H^2(X,W\OO_X)$ may not be
finitely generated for surfaces on the Noether lines
in arbitrary large characteristics, let
us recall the following result \cite[Theorem 5.4]{lie sch}:

\begin{Theorem}[--,Sch\"utt]
 There exists an arithmetic progression $P$ of primes
 of density at least $0.99999985$ such that for all $p\in P$
 there exists an even Horikawa surface $X_p$ in 
 characteristic $p$ that is unirational.\qed
\end{Theorem}

For such unirational surfaces 
$\fBr(X)$ is unipotent \cite[Remarque II.5.13]{illusie}.
Thus:

\begin{Corollary}
 For all $p\in P$ there exists an even Horikawa surface
 $X_p$ in characteristic $p$, such that
 $H^2(X_p,W\OO_{X_p})$ is not finitely generated.
 In particular, its slope spectral sequence 
 does not degenerate at $E_1$-level.
\end{Corollary}

\section{Lifting and topological flatness of moduli spaces}
\label{sec:lifting}

In this section we show that the canonical model of a surface
on the Noether lines lifts over the Witt ring $W(k)$.
Its minimal model lifts via an algebraic space
over possibly ramified extensions of $W(k)$.
Finally, we use these lifting results to prove that
moduli spaces of surfaces on the Noether lines
are topologically flat over the integers.

We start with a more general observation:
let $k$ be a field of characteristic $p>0$.
We choose a DVR
$(R,\idealm)$ of characteristic zero
with residue field $R/\idealm\iso k$ and denote its
quotient field by $K$.
Let us also assume that $R$ is a Nagata ring.
By a {\em lifting} of a variety $X\to\Spec k$ over $R$ 
we mean a scheme (or an algebraic space if explicitly stated)
%${\cal X}\to\Spec R$ 
that is flat over $\Spec R$ with special fiber $X$.

\begin{Proposition}
 \label{canonicallifting}
 Let $X\to\Spec k$ be a surface of general type, 
 not necessarily minimal and with at worst Du~Val singularities.
 If $X$ lifts over $R$ then
 \begin{enumerate}
  \item its canonical model $\Xcan$ also lifts over $R$, and
  \item there exists a possibly ramified extension $R'\supseteq R$
    and an algebraic space, flat over $R'$, with special
    fiber the minimal smooth model
    $X_{\rm min}$ of $X$.
 \end{enumerate}
\end{Proposition}

\prf
Let ${\cal X}\to\Spec R$ be a lifting of $X$.
Denote by $K$ the field of fractions of $R$ and by
${\cal X}_K$ the generic fiber of $\cal X$.
Since $X$ and $R$ are Gorenstein, the same is true 
for $\cal X$ and ${\cal X}_K$.
We denote by $\omega_{\cal X}$ its dualizing sheaf and
note its restriction to the generic and special
fiber yields the respective dualizing sheaves.
Note that ${\cal X}_K$ also has at worst 
Du~Val singularities \cite[Proposition 6.1]{lie2}.

Since $\cal X$ is flat over $R$ we have
$\chi:=\chi(\OO_X)=\chi(\OO_{{\cal X}_K})$ and
$K^2:=K_X^2=K_{{\cal X}_K}^2$.
Moreover, we have $h^1(X,\omega_X^{\otimes n})=0$ 
for $n\geq3$ by \cite[Theorem II.1.7]{ek} 
and thus Riemann--Roch yields
$$
 h^0(X,\omega_X^{\otimes n}) \,=\, \frac{n(n-1)}{2}K^2\,+\,\chi
 \,=\,
 h^0({\cal X}_K,\omega_{{\cal X}_K}^{\otimes n})
 \,\mbox{ \quad for all \quad }\,
 n\geq3.
$$
This implies that for $n\geq3$ global sections
of $\omega_X^{\otimes n}$ lift to global sections of
$\omega_{\cal X}^{\otimes n}$ and that the
$H^0({\cal X},\omega_{\cal X}^{\otimes n})$
are torsion-free $R$-modules
\cite[Corollary III.12.9]{hart}.
Thus, the $n$-fold Veronese map for some $n\geq3$ 
of the relative canonical ring
$$
{\cal X}_{\rm can}^{(n)}
\,:=\,
\Proj\,\bigoplus_{k=0}^\infty\, H^0({\cal X},\,\omega_{\cal X}^{\otimes nk}) 
\,\to\,\Spec R
$$
is a projective scheme over $R$ with special (resp. generic) 
fiber isomorphic to the canonical model of $X$ (resp. ${\cal X}_K$). 
This gives the stated lifting of $\Xcan$.

By the main result of \cite{Artin} there exists a possibly
ramified extension $R'$ of $R$ and an algebraic space over
$\Spec R'$ over which ${\cal X}_{\rm can}^{(n)}\to\Spec R$
admits a simultaneous and minimal resolution of singularities.
This gives the desired lifting of $X_{\rm min}$. 
\qed\medskip

Using the explicit classification results of the previous sections
we obtain

\begin{Theorem}
 Let $X$ be a surface on the Noether lines over 
 an algebraically closed field $k$ of positive characteristic.
 Then
 \begin{enumerate}
  \item the canonical model $\Xcan$ lifts over $W(k)$ 
  \item the minimal smooth model $X_{\rm min}$ lifts via an algebraic
    space over a possibly ramified extension of $W(k)$.
 \end{enumerate}
\end{Theorem}

\prf
We only deal with the case that $X$ is an even Horikawa
surface or an odd Horikawa surface with $p_g\geq5$
and leave the remaining cases to the reader.

From Theorem \ref{thm:evenstructure} and 
Theorem \ref{thm:oddstructure}
we see that there exists a surface $X'$, birational to $X$, 
with at worst Du~Val singularities and a finite flat
morphism $\pi:X'\to S$, where $S$ is a smooth rational
surface. 
As usual, we denote by ${\cal L}$ the invertible sheaf
associated to $\pi$.
From the vanishing 
${\rm ext}^1({\cal L}^\vee,\OO_S)=h^1({\cal L})=0$
(Lemma \ref{lemma:vanishing})
we see that $X'\iso\uSpec(\OO_S\oplus{\cal L}^\vee)$
as $S$-scheme.

As already discussed in \cite[Section 7]{lie2},
Grothendieck's existence theorem implies that rational
surfaces, and in particular $S$, lift over $W(k)$.
Moreover, invertible sheaves on rational surfaces
even lift uniquely.
We can lift sections of these invertible sheaves if
their $h^1$ vanishes \cite[Lemma 7.1]{lie2}.
This being the case for 
${\cal L}$ and ${\cal L}^{\otimes 2}$ 
by Lemma \ref{lemma:vanishing},
it follows that the whole double cover 
$\pi:X'=\uSpec(\OO_S\oplus{\cal L}^\vee)\to S$ 
lifts over $W(k)$.

The result now follows from Proposition \ref{canonicallifting}.
\qed\medskip

Let us recall that a scheme over a regular, integral, one-dimensional
base scheme is called {\em topologically flat} 
if its generic fiber is dense. 

\begin{Theorem}
  \label{thm:flat}
  The moduli spaces of surfaces on the Noether lines are topologically
  flat over $\Spec\ZZ$.
\end{Theorem}

\prf
This follows immediately from the lifting results above.
\qed\medskip

\begin{Remark}
 Of course, it would be nicer to prove flatness of these moduli spaces
 and to understand their scheme structure (singularities, reducedness)
 in detail.
 This is equivalent to understanding the deformation theory of our
 surfaces, e.g., as in \cite{horquint}.
 Unfortunately, it seems that the contribution of Du~Val singularities 
 on their canonical models to their deformation 
 functors is more complicated than in characteristic zero.
 We will come back to this in another article.
\end{Remark}

\begin{appendix}

\section{Double covers and canonical resolution of singularities}
\label{sec:appendix}

In this appendix we recall a couple of facts  
on double covers and canonical resolution of their singularities.
These results are somewhat scattered over the literature and divided
into several subcases, depending on whether the characteristic is
different from $2$ or not.
In characteristic $2$, the setup is usually divided
into the separable and the inseparable case.
Here, we provide the material needed to make the arguments
of \cite[Section 2]{horquint} and \cite[Section 1]{horikawa2}
work in arbitrary characteristic.

First, a characteristic-free description of double covers 
is given in \cite[Chapter 0.1]{cd}.
In characteristic $\neq2$, double covers and resolution of their singularities 
can be described in terms of branch loci,
cf. \cite[Section 2]{horquint} and \cite[Section III.7]{bhpv}.
In characteristic $2$, there are two distinct types of
double covers:
they are either separable in which case we are dealing with Artin--Schreier 
extensions (of simple type), cf. \cite{tak}, or they are
purely inseparable and we are dealing with $p$-closed foliations 
and their singularities, cf. \cite{hir}.
These latter two cases can be uniformly described in the language
of $\alpha_s$-torsors, see \cite[I.(1.10)]{ek} or
Section \ref{sec:alphaltorsors}.

Let $f:X\to S$ be a generically finite morphism of degree $2$ from
a normal onto a smooth surface.
Let $\nu:X\to X'$ followed by $f':X'\to S$ be
the Stein factorization of $f$.
Since $X'$ is normal and $f'$ is finite, it follows
that $f'$ is flat.
Moreover, since $f'$ is of degree $2$ we obtain a short exact sequence
\begin{equation}
 \label{associated}
 0\,\to\,\OO_S\,\to\,f'_\ast\OO_{X'}\,\to\,{\cal L}^\vee\,\to\,0
\end{equation}
where ${\cal L}^\vee$ is an invertible sheaf on $S$.
We define $\cal L$ to be its dual.

\begin{Definition}
 We refer to $\cal L$ as the \emph{invertible sheaf associated with $\pi$}.
\end{Definition}

In characteristic $\neq2$, the morphism $f'$ is Galois and the decomposition
of $f'_\ast\OO_{X'}$ into eigensheaves yields a splitting of
(\ref{associated}).
Moreover, the branch locus of $f'$ corresponds to a section of
${\cal L}^{\otimes2}$ and defines the $\OO_S$-algebra structure on 
$f'_\ast\OO_{X'}$ uniquely.
However, in characteristic $2$, there is no notion
of branch locus if $f'$ is inseparable.
In the separable case this branch locus corresponds to a section
of $\cal L$, which is in general not enough to determine  
the $\OO_S$-algebra structure of $f'_\ast\OO_{X'}$.

Since $f'$ is flat and its fibers are Artin-algebras of length $2$,
which are Gorenstein, it follows that $X'$ is Gorenstein.
More precisely, we have
\begin{equation}
 \label{canonical}
 \omega_{X'}\,\iso\,f'^\ast(\omega_S\otimes{\cal L}),
 \mbox{ \quad hence \quad }
 K_{X'}^2\,=\,\frac{1}{2}(K_S+L)^2\,,
\end{equation}
where $\omega_{X'}$ denotes the dualizing sheaf of $X'$,
cf. \cite[Section 0.1]{cd}.
Since $X'$ is Gorenstein, there exists a Cartier divisor on $X'$,
the {\em canonical cycle} $Z_K$, which is
supported on the exceptional locus of $\nu$ and such that
$$
\omega_X\,\iso\,\nu^\ast\omega_{X'}\otimes\OO_X(-Z_K)\,.
$$
If $\nu$ is {\em minimal}, i.e., 
does not contract any $(-1)$-curves, then
$Z_K\geq0$.
%, cf. \cite[Section 4.19]{reid} or \cite[Lemma 5]{horquint}.
Then,
$Z_K$ being zero is equivalent to $X'$ having only rational singularities,
i.e., Du~Val singularities in our setup.

We already remarked that ${\cal L}^{\otimes2}$ replaces the notion of
a branch divisor.
Then we define an invertible sheaf $\cal R$ by
$$
 \omega_X \,\iso\, f^\ast\omega_S\otimes{\cal R}
 \mbox{ \quad such that the familiar formula \quad }
 \OO_X(Z_K)\,\iso\,f^\ast{\cal L}\otimes{\cal R}^\vee
$$
from \cite[Lemma 4]{horquint} holds true.
This $\cal R$ plays the role of a ramification divisor, 
although there is in general no distinguished section
of this invertible sheaf.

To resolve the singularities of a finite and flat 
double cover, there is the following procedure, which
is classical in characteristic $\neq2$, cf. 
\cite[Section III.7]{bhpv} and has been shown in 
characteristic $2$ in the separable case in
\cite[Section 2]{tak} and in the inseparable case in
\cite[Proposition 2.6]{hir}.

\begin{Proposition}[canonical resolution of singularities]
 \label{prop:canres}
 Let $f:X\to S$ be a finite and flat morphism of degree $2$
 from a normal onto a smooth surface.
 Let $x\in X$ be a singular point.
 Then there exists a sequence of blow-ups in smooth points
 $\pi:\tilde{S}\to S$ such that the normalized fiber product
 $$\begin{array}{lcl}
    X&\stackrel{\varpi}{\leftarrow}&\widetilde{X}\\
    \downarrow {\scriptstyle{f}}&&\downarrow{\scriptstyle{{\widetilde{f}}}}\\
    S&\stackrel{\pi}{\leftarrow}&\widetilde{S}
   \end{array}$$
 yields a resolution of singularities $\varpi:\widetilde{X}\to X$
 and such that $\widetilde{f}$ is again a finite and flat morphism 
 of degree $2$.\qed
% In particular, a resolution of singularities of $X$ can be achieved within
% the category of finite and flat double covers.
\end{Proposition}

We need a characteristic-free substitute for 
\cite[Lemma 6]{horquint}.
Thus, we want to determine the effect of
a single blow-up and normalized fiber product in the
previous proposition, which in general is only
a partial desingularization of $X$.

\begin{Lemma}
 \label{invariants}
 In the situation of Proposition \ref{prop:canres}, assume that
 $\pi$ is only a single blow-up in $f(x)$.
%, i.e.,
% $\varpi:\widetilde{X}\to X$ is possibly only
% a partial desingularization.
 Let $E$ be the exceptional divisor of $\pi$ and
 $\widetilde{{\cal L}}$ the invertible sheaf
 associated with $\tilde{f}$.
 Then
 there exists an integer $n\geq0$ such that
 $$
 \begin{array}{ccl}
   \blowup\widetilde{{\cal L}} &\iso& \pi^\ast{{\cal L}}\otimes\OO_{\widetilde{S}}(-nE)\\
   \blowup\omega_{\widetilde{X}} &\iso& \varpi^\ast(\omega_X)\otimes \widetilde{f}^\ast\OO_{\widetilde{S}}((1-n)E)\\
   \blowup K_{\widetilde{X}}^2 &=& K_X^2 \,-\, 2(n-1)^2\\
   \blowup \chi(\OO_{\widetilde{X}}) &=& \chi(\OO_X) \,-\, \frac{1}{2}n(n-1)\\
   \blowup &=&2\chi(\OO_S)\,+\,\frac{1}{2} {\cal L}\cdot (\omega_S\otimes{\cal L}) \,-\,\frac{1}{2}n(n-1)\,. 
 \end{array}
 $$
 Moreover,
 \begin{enumerate}
  \item $x\in X$ is a smooth point if and only if $n=0$, and
  \item if $x\in X$ is a rational singularity then $n=1$.
 \end{enumerate}
\end{Lemma}

\prf
Pulling back (\ref{associated}) to $\widetilde{S}$, we obtain a commutative diagram 
\begin{equation}
 \label{lemmadiagram}
 \begin{array}{ccccccccc}
  0&\to&\pi^\ast\OO_S&\to&\pi^\ast f_\ast\OO_X&\to&\pi^\ast{\cal L}^\vee&\to&0\\
  &&\downarrow{\scriptstyle{\iso}}&&\downarrow&&\downarrow\\
  0&\to&\OO_{\widetilde{S}}&\to&\widetilde{f}_\ast\OO_{\widetilde{X}}&\to&\widetilde{\cal L}^\vee&\to&0
 \end{array}
\end{equation}
where $\pi^\ast f_\ast\OO_X \to \widetilde{f}_\ast\OO_{\widetilde{X}}$ is the normalization.
The induced inclusion $\pi^\ast{\cal{L}}^\vee\to\widetilde{{\cal L}}^\vee$ is an 
isomorphism outside $E$. 
Hence there exists an integer $n\geq0$ such that
$$
\pi^\ast {\cal L}^\vee\,\iso\,\widetilde{\cal L}^\vee\otimes\OO_{\widetilde{S}}(-nE)\,.
$$
The formulae for $\omega_{\widetilde{X}}$ and $K_{\widetilde{X}}^2$ follow 
from (\ref{canonical}).
We compute $\chi(\OO_{\widetilde{X}})$ by
taking Euler characteristics in (\ref{lemmadiagram}) and applying Riemann--Roch.

If $x\in X$ is singular, then $X\times_S\widetilde{S}$ is singular
along $\varpi^{-1}(x)$.
In particular, $\pi^\ast f_\ast\OO_X$ is not normal and we get $n\geq1$.

If $x\in X$ has at worst rational singularities then so has $\widetilde{X}$
along $\varpi^{-1}(x)$ and thus $\chi(\OO_X)=\chi(\OO_{\widetilde{X}})$, 
which gives $n\in\{0,1\}$.

Finally, if $x\in X$ is a smooth point, then $\widetilde{X}$ is smooth along
$\varpi^{-1}(x)$, which implies that $K_{\widetilde{X}}^2<K_X^2$.
In particular, $n\neq1$ which gives $n=0$ be the previous consideration.
\qed\medskip

\begin{Remark}
  In characteristic $\neq2$, the integer $n$ is half the 
  multiplicity of a singularity of the branch curve,
  see \cite[Lemma 6]{horquint} or \cite[Theorem III.(7.2)]{bhpv}.
\end{Remark}

We note that our singularities are automatically Gorenstein.
Thus, rational singularities on double covers are
Du~Val singularities.
The following extends \cite[Lemma 5]{horquint}
and \cite[Theorem III.(7.2)]{bhpv} 
to arbitrary characteristic: 

\begin{Proposition}
 Let $f:X\to S$ be a finite and flat morphism of degree $2$,
 where $S$ is a smooth surface and $X$ has a 
 Du~Val singularity in $x\in X$.
 Then the canonical resolution 
 $\varpi:\widetilde{X}\to X$ coincides with the minimal 
 resolution of singularities.
 Moreover,
 $$
    \widetilde{{\cal L}} \,\iso\, \pi^*{\cal L}\otimes\omega_{\widetilde{S}}^\vee\otimes\pi^*\omega_S\,.
 $$
% and
%\marginpar{which cycle?}
% $$
%    \OO_{\widetilde{X}}(Z_{\rm num})\,\iso\,\widetilde{f}^\ast\left(
%     \omega_{\widetilde{S}}\otimes\pi^*\omega_S^\vee \right)\,,
% $$ 
% where $Z_{\rm num}$ is the numerical, or, fundamental cycle of the singularity.
% WRONG already for A_{2k+1}-singularities in characteristic zero we get strange things:
% numerical cycle would be just sum over the (-2)-curves
% the K_{\widetilde{S}}-K_S is the sum \sum_i iE_i, where the E_i are from left to
% right as in \cite[Table 1, Section III.7]{bhpv}
\end{Proposition}

\prf
By the previous lemma, we find $n=1$ for every single blow-up in a 
closed point of the base. 
Thus, the canonical resolution of singularities fulfills
$\omega_{\widetilde{X}}=\varpi^*\omega_X$, which implies that
it is minimal.

From the formulae
$$
 \omega_X\,\iso\,f^*(\omega_S\otimes{\cal L}),\,\mbox{ \quad }\,
 \omega_{\widetilde{X}}\,\iso\,\widetilde{f}^*(\omega_{\widetilde{S}}\otimes\widetilde{{\cal L}})\,
 \mbox{ \quad and \quad }\,
 \omega_{\widetilde{X}}\,\iso\,\varpi^\ast(\omega_X)
$$
we obtain the statement linking $\pi^\ast{\cal L}$ to $\widetilde{{\cal L}}$. 
\qed

\end{appendix}


\begin{thebibliography}{EGA IV}
  \bibitem[Ar]{Artin} M.~Artin, {\it Algebraic Construction of Brieskorn's 
    Resolutions}, J. of Algebra 29, 330-348 (1974).
  \bibitem[AM]{artin mazur} M.~Artin, B.~Mazur, {\em Formal groups arising from 
    algebraic varieties}, Ann. Sci. \'Ecole Norm. Sup. 10, 87-131 (1977).
  \bibitem[AK]{ashikaga konno} T.~Ashikaga, K.~Konno, {\it Algebraic surfaces of 
    general type with $c^2_1=3p_g-7$}, Tohoku Math. J. 42, 517-536 (1990).
  \bibitem[BHPV]{bhpv} W.~P.~Barth, K.~Hulek, C.~Peters, A. van de Ven,
     {\it Compact Complex Surfaces}, 2nd edition, Erg. d. Math., 3. Folge
     Volume 4, Springer (2004). 
  \bibitem[Bea]{bea} A.~Beauville, {\it L'application canonique pour les surfaces
     de type g\'en\'eral}, Invent. Math. 55, 121-140 (1979).
  \bibitem[Bo]{bom} E.~Bombieri, {\it Canonical models of surfaces of general
     type}, Inst. Hautes \'Etudes Sci. Publ. Math. 42, 171-220 (1973).
  \bibitem[C]{castelnuovo} G.~Castelnuovo, {\it Osservazioni intorno alia geometria 
     sopra una superficie, Nota II}, Rendiconti del R. Instituto Lombardo 24 (1891).
  \bibitem[CD]{cd} F.R.~Cossec, I.V.~Dolgachev, {\it Enriques Surfaces I},
     Prog. in Math. 76, Birkh\"auser (1989).
  \bibitem[DI]{deligne illusie} P.~Deligne, L.~Illusie, 
     {\it Rel\`evements modulo $p^2$ et d\'ecomposition du complexe de de Rham},
     Invent. Math. 89, 247-270 (1987).
  \bibitem[EH]{eh} D.~Eisenbud, J.~Harris, {\it On Varieties of Minimal Degree
     (A Centennial Account)}, Algebraic Geometry, Bowdoin 1985, Proc. Symp. Pure
     Math. 46, Part 1, 3-13 (1987).
%  \bibitem[EKS]{eks} D.~Eisenbud, J.~Koh, M.~Stillman, {\it Determinantal 
%     equations for curves of high degree},
%     Am. J. Math. 110, 513-539 (1988).
  \bibitem[Ek1]{ek gauge} T.~Ekedahl, {\em Diagonal complexes and $F$-gauge structures},
     Hermann  (1986).
%  \bibitem[E]{ek fol} T.~Ekedahl, {\it Foliations and Inseparable Morphisms}, 
%     Proc. Symp. Pure Math. 46, No. 2, 139-149 (1987).
  \bibitem[Ek2]{ek} T.~Ekedahl, {\it Canonical models of surfaces of general type in 
     positive characteristic},  Inst. Hautes \'Etudes Sci. Publ. Math.  No. 67, 
     97-144 (1988).
%  \bibitem[El]{el} R.~Elkik, {\it Singularit\'es rationelles et d\'eformations},
%     Invent. math. 47, 139-147 (1978).
  \bibitem[En]{enr} F.~Enriques, {\it Le superficie algebriche},
     Nicola Zanichelli (1949).
%  \bibitem[EGA IV]{ega4} A.~Grothendieck, {\it \'El\'ements de g\'eom\'etrie 
%     alg\'ebrique IV: \'Etude locale des sch\'emas et des morphismes de
%     sch\'emas}, Publ. Math. IHES 20 (1964), 24 (1965), 28 (1966), 32 (1967).
  \bibitem[Hart]{hart} R.~Hartshorne, {\it Algebraic Geometry}, Springer (1977).
  \bibitem[Hir]{hir} M.~Hirokado, {\it Singularities of multiplicative $p$-closed
     vector fields and global $1$-forms on Zariski surfaces}, J.~Math.~Kyoto Univ.,
     39-3, 455-468 (1999)
% \bibitem[Hir2]{hir} M.~Hirokado, {\it Zariski surfaces as quotients of Hirzebruch
%     surfaces by 1-foliations}, Yokohama Math. J. 47, 103-120 (2000).
  \bibitem[HorQ]{horquint} E.~Horikawa, {\it On Deformations of Quintic Surfaces},
     Invent.~math. 31, 43-85 (1975).
  \bibitem[Hor1]{horikawa} E.~Horikawa, {\it Algebraic surfaces of general
     type with small $c_1^2$, I}, Ann. Math. 104, 357-387 (1976).
  \bibitem[Hor2]{horikawa2} E.~Horikawa, {\it Algebraic surfaces of general 
     type with small $c_1^2$, II}, Invent. Math. 37, 121-155 (1976).
  \bibitem[Ill]{illusie} L.~Illusie, {\it  Complexe de de~Rham--Witt et 
     cohomologie cristalline}, Ann. Sci. \'Ecole Norm. Sup. 12, 501-661 (1979).
  \bibitem[Ill2]{ill2} L.~Illusie, {\it Grothendieck's existence theorem in formal
    geometry with a letter of Jean-Pierre Serre}, 
    AMS Math. Surveys Monogr. 123, Fundamental Algebraic Geometry, 
    179-233 (2005). 
  \bibitem[Jou]{jou} J.-P.~Jouanolou, {\it Th\'eor\`emes de Bertini et
     Applications}, Prog. in Math. 42, Birkh\"auser (1983).
  \bibitem[Kar]{karras} U.~Karras, {\em Local cohomology along exceptional sets},
     Math. Ann. 275, 673-682 (1986).
  \bibitem[Ka]{Katz} N.~Katz, {\it Nilpotent connections and the monodromy theorem:
     applications of a result of Turrittin}, Publ. Math. IHES 39, 175-232 (1970).
  \bibitem[Lie1]{lie} C.~Liedtke, {\it Uniruled surfaces of general type},
      Math. Z. 259, 775-797 (2008).
  \bibitem[Lie2]{lie2} C.~Liedtke, {\it Algebraic surfaces of general type with 
     small $c_1^2$ in positive characteristic},
     Nagoya Math. J. 191, 111-134 (2008).
  \bibitem[Lie3]{lie godeaux} C.~Liedtke, {\it Non-classical Godeaux Surfaces},
     Math. Ann. 343, 623-637 (2009).
  \bibitem[LS]{lie sch} C.~Liedtke, M.~Sch\"utt, {\it Unirational Surfaces on the 
     Noether Line}, Pacific J. Math. 239, 343-356 (2009). 
  \bibitem[Lip]{lip} J.~Lipman, {\it Rational singularities, with applications
     to algebraic surfaces and unique factorization}, 
     Publ. Math. IHES 36, 195-279 (1969).
  \bibitem[Lip2]{lip2} J.~Lipman, {\it Desingularization of two-dimensional 
     schemes}, Ann. Math. 107, 151-207 (1978).
  \bibitem[Mi]{mir} R.~Miranda, {\it Nonclassical Godeaux surfaces in 
    characteristic five}, Proc. Amer. Math. Soc.  91, 9-11  (1984).
  \bibitem[Noe]{noether} M.~Noether, {\it Zur Theorie der eindeutigen Entsprechungen
      algebraischer Gebilde}, Math. Ann. 8, 495-533 (1875).
   \bibitem[Pa]{pardini} R.~Pardini, {\it Triple covers in positive characteristic},
       Ark. Mat. 27, 319-341 (1989).
   \bibitem[Ram]{ram} C.P.~Ramanujam, {\it Remarks on the Kodaira vanishing theorem},
      Journal of the Indian Math. Soc. 36, 41-51 (1972).
%   \bibitem[Rei]{reid} M.~Reid, {\it Chapters on Algebraic Surfaces},
%      Complex algebraic geometry,
%      IAS/Park City Math. Ser., 3, Amer. Math. Soc., 3-159 (1997).
%   \bibitem[R-S]{rs} A.~N.~Rudakov, I.~R.~\v{S}afarevi\v{c}, {\em 
%      Inseparable Morphisms of Algebraic Surfaces}, 
%      Izv. Akad. Nauk SSSR Ser. Mat. 40, 1269-1307 (1976).
%   \bibitem[S]{saf} I.~R.~\v{S}afarevi\v{c} et al., {\it Algebraic Surfaces},
%      Proc. Steklov Inst. 75 (1965).
   \bibitem[S-B]{sb} N.I.~Shepherd-Barron, {\it Unstable vector bundles and linear
      systems on surfaces in characteristic $p$}, Invent. math. 106, 243-262 (1991).
%   \bibitem[S-B2]{sb2} N.I.~Shepherd-Barron, {\it Some foliations on surfaces in
%      characteristic $2$}, J. Alg. Geom. 5, 521-535 (1996).
   \bibitem[Ta]{tak} Y.~Takeda, {\it Artin--Schreier coverings of algebraic surfaces},
      J.~Math.~Soc. Japan 41, 415-435 (1989).
\end{thebibliography}
\end{document}